\documentclass[a4paper,notitlepage, 12pt,reqno]{article}
  \usepackage{graphicx}
  \usepackage{mathtext}
 \usepackage[cp1251]{inputenc}
  \usepackage[T2A]{fontenc}
 \usepackage[russian,english]{babel}
 \usepackage[all,arc,poly,2cell,curve,arrow,tips]{xy}
  \usepackage{enumerate, eucal, amsthm,amsmath, amssymb}

  \allowdisplaybreaks[4]

 \theoremstyle{definition}
 \newtheorem{defn}{Definition}%[section]

 \theoremstyle{plain}
 \newtheorem{thm}{Theorem}
 \newtheorem*{thm*}{Theorem}
 \newtheorem{prop}{Proposition}
  \newtheorem*{prop*}{Proposition}
 \newtheorem{cor}{Следсвие}
  \newtheorem*{cor*}{Следствие}
 \newtheorem{lem}{Lemma}
  \newtheorem*{lem*}{Лемма}

 \theoremstyle{remark}

 \newtheorem*{remark*}{Замечание}

 \renewcommand{\abstractname}{}

  \newcounter{ab}
\title{Wigner coefficient for  Lie algebras of series $B,C,D$ and a base of Gelfand-Tsetlin type.\thanks
 {The work was supported by grants  NSch-5998.2012.1, RFFI-12-01-31414, MK-4594.2013.1} \\
 \begin{small}
Dedicated to the 100-th anniversary of I.M. Gelfand
\end{small}
}

 \author{D.V. Artamonov\footnote{artamonov.dmitri@gmail.com}, V.A. Golubeva\footnote{goloubeva@yahoo.com}}
  \date{}

\begin{document}
 \maketitle

\renewcommand{\abstractname}{}

\begin{abstract} For the Lie algebras $g_n= \mathfrak{o}_{2n+1},\mathfrak{sp}_{2n},\mathfrak{o}_{2n}$
a simple construction of a base in an irreducible representation is
given. The construction of this base uses the method of
$Z$-invariants of Zhelobenko and the technique of Wigner
coefficients, which was applied by Biedenharn and Baird to the
construction of a Gelfand-Tsetlin base in the case
$\mathfrak{gl}_n$. A relation between matrix elements and Wigner
coefficients for
  $g_n$ and analogous objects for $\mathfrak{gl}_{n+1}$ is established.

\end{abstract}

%\tableofcontents

\section{Introduction}

One can construct explicitly a representation of a  simple Lie
algebra of the type $A$  using the technique of Wigner coefficients.
In the case $\mathfrak{sl}_2$  such a construction was obtained by
 Wigner and Racah, and in the case
$\mathfrak{sl}_n$ it was obtained by Biedenharn, Baird, Louck and others.  The aim
of the present paper is a generalization of this technique to the
case of simple Lie algebras $g_n $ of types $B,C,D$.  Such
constructions
 can be applied in the nuclear physics, for example in the theory of nuclear shells.

The construction given in the present paper establishes a remarkable relation between Gelfand-Tsetlin bases for the algebras of series
 $A$ and of series $B,C,D$. Moreover we express matrix elements for generators in the base and Wigner coefficients for algebras of series  $B,C,D$ through the analogous objects of series $A$.

The construction in the paper is inductive. For the algebras
$g_n=\mathfrak{o}_{2n+1},\mathfrak{sp}_{2n},\mathfrak{o}_{2n}$ the
constructions are different in the case $n=1$ and also in the case
$n=2$. But the consideration  for $n>2$  are similar for all algebras
and also they are similar to  analogous constructions for the
algebra $\mathfrak{gl}_{n+1}$. This fact plays a key role in all
constructions. It is a direct corollary of the fact that for the
Dynkin diagrams $A_{n+1},B_{n},C_{n},D_{n}$   transform in a similar
way when we change  $2$ to $n$.
% (see figure \ref{fig1}).

%\begin{figure}
%\centering
%\includegraphics[bbllx=0,bblly=0,bburx=120,bbury=100]{dd.bmp}
%\caption{The Dynkin diagramms.}\label{fig1}
%\end{figure}

The construction uses a well-known inductive procedure of
Gelfand and Tsetlin. The main step in this construction is an
investigation of a branching of an irreducible representation of
$g_n$ when one restricts the algebra  $g_{n}\downarrow g_{n-1}$. For this
investigation we use the method of $Z$-invariants of Zhelobenko.

To obtain formulas for the action of generators of the algebra in
the base we use the technique of Wigner coefficients. More precise, the matrix
elements are expressed through the simplest Wigner coefficients that
correspond to a decomposition of a  tensor product of an arbitrary
representation with a standard representation of the algebra.
Earlier analogous results were obtained by Biedenharn and Baird.
They   expressed  matrix elements for operators $E_{i,i-1}$ in the
Gelfant-Tsetlin base for the algebras of series $A$ through Wigner coefficients \cite{1963}.

I.M. Gelfand in the late 80-th several times said about the desire
to generalize the works of
 Biedenharn and Louck. We have understood this as a problem of a generalization of the construction
  of Biedenharn and Louck for the  Lie algebras of series $A$ to the case of  algebras of series
   $B,C,D$. This is done in the present paper.
%    \footnote{Note that similar questions were considered in the papers of Gelfand and Tsentlin but they constructed  not a base for series  $B,C,D$ but a base whose construction uses restrictions $\mathfrak{o}_N\downarrow \mathfrak{o}_{N-1}$.}.

The problem of construction of  base for orthogonal algebras
algebras   was investigated in the works of Gelfand
and Tsetlin, but  it's construction is based on  restrictions $\mathfrak{o}_N\downarrow
\mathfrak{o}_{N-1})$ \footnote{Thus this base is not a base for a series $B$ or $D$ since the algebras of both series are involved in it's construction}.
%they obtained branching rulers for irreducible
%representations, when one restricts $g_n\downarrow g_{n-1}$.  This
%gives an indexation of base vectors. But for a long time one could
%not obtain formulas for the action of generators.

In \cite{zh} using the method  of $Z$-invariants of
Zhelobenko  there was constructed a Gelfanf-Tsetlin type base for
orthogonal algebras (based on and for symplectic algebras
 (see also
 \cite{Sch1}, \cite{Sch2}, \cite{Sch3}).  But in these papers only an indexation is constructed, the formulas for the action of generators are not obtained.

Nevertheless there exists a construction of a Gelfand-Tsetlin type
base for the algebras of series $B,C,D$ which belongs to Molev
\cite{M}.  But it uses a difficult technique of
Mickelsson-Zhelobenko algebras and the action of Jangians on the
multiplicity spaces\footnote{ An indexation of base vectors in the base that is constructed in the present paper is the same as
 in the Gelfand-Tsetlin-Molev \cite{M}. But we have not managed to construct an exlicit isomorphism betweeen
 our base and the Gelfand-Tsetlin-Molev base}.

The construction given in the present paper is
much simpler.  Is's main advantage is that it establishes a remarkable  relation between Gelfant-Tsetlin bases for
series $A,B,C,D$. Also as a by-product of our construction we obtain
explicit formulas for Wigner coefficients for algebras of series
$B,C,D$. The indexation of base vectors is the same as Molev's.  But
also our base is not orthonormal as Molev's.

The  Gelfand-Tsetlin type bases for Lie $g_n $ of
types $B_n,C_n,D_n$ based on restrictions $g_n\downarrow g_{n-1}$
are important in nuclear physics.  Such a base is used in
problems where the algebra of five-dimensional  quasi-spin (isomorphic to
$\mathfrak{o}_5=\mathfrak{sp}_4$) is involved, for example in the
theory of nuclear shells, \cite{He}, in the Bohr-Mottelson model
\cite{BM}, in the model of interacting bosons \cite{IBM}, in the
models of high-temperature superconductivity \cite{16}, \cite{17},
\cite{GoLi}.

 The typical application of Wigner coefficients is the following. If one
  identify an irreducible representation with a (quasi)particle, then Wigner
  coefficients that define a decomposition of a tensor product of two representation
  describe a spectrum of (quasi)particles that appear  after their  interaction
(см. \cite{Lip}).

The Wigner coefficients that are calculated in the present paper correspond to the adding of one (quasi)particle. Thus in the case of a nuclear shell model the multiplication to the standard representation corresponds to the adding to the system of one nucleon that is a proton or a neutron. The explicit formulas for Wigner coefficients allow to write the selection rulers for the quantum numbers that define  the states of quasi-spin and calculate the probabilities of transitions into these states.

There exist also completely different applications. Thus the
high-temperature decomposition in the classical  $N$-vector model is
described using Wigner coefficients of the orthogonal algebra
$\mathfrak{o}_N$ \cite{joyce}.

In all these problems mostly the Wigner coefficients that define the
decomposition of two symmetric representations are used. Such
coefficients were explicitly obtained in many particular cases in
\cite{Gav}, \cite{Ki}, \cite{Junkr}, \cite{Al1}, \cite{Al2}. But mention
the authors of these papers used the Gelfand-Tsetlin base whose
construction  is based on restrictions $\mathfrak{o}_N\downarrow\mathfrak{o}_{N-1}$
 ($\mathfrak{o}_5\downarrow\mathfrak{o}_4\downarrow\mathfrak{o}_3$ in the case $\mathfrak{o}_5$).
This base is not natural form the physical point of view. In the
present paper we use the base for a representation of
$\mathfrak{o}_5$, whose construction is based on restrictions within
the series $B$.

\subsection{The structure of the paper}

The structure of the paper is the following. In Section \ref{s2}
definitions and notations are introduced.

The technical details, the discussion of definitions is placed in
Appendix \ref{appendi}.

In Section \ref{oldvector} the method of $Z$-invariants is
explained, it allows to solve effectively the problem of description
of a branching of a  representation when one restricts an algebra.

In Section \ref{shema} we give the scheme of a solution of the problem
of restriction $g_n\downarrow g_{n-1}$\footnote{When we write "the
problem of restriction $g_n\downarrow g_{n-1}$"  we mean the problem
of an explicit description of a base in the space of
$g_{n-1}$-highest vectors in a
 $g_{n}$-representation.}.  We explain how to construct a Gelfand-Tselin type
 base, obtain coefficients and restricted  Wigner coefficients and
 matrix elements of the generators.

The main construction is given in Sections \ref{step1}, \ref{step2},
\ref{step3}. In these Section the Gelfand-Tsetlin type base for
series $B,C,D$ and explicit formulas for the action of generators
are constructed for  $n=1$, $n=2$, $n>2$ respectively.  Note that
Wigner coefficients for $g_{n-1}$  are calculated when the algebra
$g_n$ is considered.

The  cases $n=1$, $n=2$  are considered in Sections
 \ref{step1}, \ref{step2}  separately for each series $B,C,D$.
 In Section \ref{step3} further steps are discussed, they are similar for all series.

In all cases all values are expressed through the analogous values
for the algebras $\mathfrak{gl}_N$,  which were obtained in an explicit
form by Biedenharn and Baird in \cite{1963} and also by Zhelobenko
in \cite{zh}.

In  Appendix the facts from the representation theory that are not
well-known are given. Mostly that can be found in \cite{1968}. The
reader can find them here of in Appendix  \ref{appendi}.

\section{Basic definitions}

\label{s2}

In the present Section the basic definitions and notations are introduced. See also Appendix  \ref{appendi}.

In the paper we use the Lie algebras $\mathfrak{sp}_{2n}$ and $\mathfrak{o}_{N}$
in split realization. These algebras act in the space with coordinates $x_{-n},...,x_{-1},x_{1},...,x_{n}$, in the cases  $\mathfrak{sp}_{2n}$ and $\mathfrak{o}_{2n}$, and in the space
 with coordinates  $x_{-n},...,x_{-1}, x_{0},x_{1},...,x_{n}$ in the case  $\mathfrak{o}_{2n+1}$.
 The generators of the symplectic algebra are the matrices $$F_{i,j}=E_{i,j}-sign(i)sing(j)E_{-j,-i},$$
 and the generators of the orthogonal algebra are the matrices  $$F_{i,j}=E_{i,j}-E_{-j,-i}.$$

\subsection{Gelfand-Tsetlin tableaux}
\label{baspred}
%Пусть $g_n$ есть либо $\mathfrak{o}_{2n+1}$, либо $\mathfrak{o}_{2n}$, либо $\mathfrak{sp}_{2n}$.

In the presen paper a Gelfand-Tsetlin base in an irreducible representation of
$g_n$ is constructed. Base vectors are indexed by Gelfand-Tsetlin  type tableaux,  these tableaux have  similar structure, let us describe it here.

% Базис
%Гельфанда - Цетлина - Молева есть базис в представлении, чья
%конcтрукция основана на ограничениях $g_n\downarrow g_{n-1}$. То
%есть, на ограничения $B_n\downarrow B_{n-1}$, $C_{n}\downarrow
%C_{n-1}$, $D_n\downarrow D_{n-1}$.

Let us be given an irreducible representation  $V$ of the algebra $g_n$ with the highest weight $[m_{-n,n},...,m_{-1,n}]$.  Then the base vectors are indexed by tableaux   $(m)$ of type

\begin{align}\label{gcmspo}
(m)=&\begin{pmatrix} [m]_{n}\\
[m']_{n}\\
[m]_{n-1}\\
[m']_{n-1}\\
...\\
[m]_1\\
[m']_{1}\\
\end{pmatrix}.
  \end{align}

The row $[m]_{n}$  is the highest weight of the considered representation
$V$,   the row $[m]_{n-1}$ is the highest weight of an irreducible
$g_{n-1}$-representation that contains the vector $(m)$. The row    $[m']_{n-1}$ is a base element in the space of  $g_{n-1}$-highest vectors with highest weight
$[m]_{n-1}$ and so on. The structure of the rows depends on the series $B_n$, $C_n$ or $D_n$.

The weight of the vector $(m)$ is denoted as
$\Delta(m)$.

%Введём ещё некоторые обозначения.

%\begin{defn}

% \end{defn}
% \begin{defn}
Denote the tableau in which all indices take the maximum possible values as $(m)_{max}$.
%\end{defn}
The vector  $(m)_{max}$ is a highest vector.

\subsection{Wigner coefficients and reduced Wigner coefficients}\label{coefvigner}

Let us define Wigner coefficients, reduced Wigner coefficients and let us give a solution of the multiplicity problem. See also Appendix  \ref{appendi}.
%Подробности могут быть найдены в приложении
%\ref{apl1}.

Denote an irreducible representation with highest weight $[m]_{n}$ as
$V^{[m]_{n}}$.
The tableau $(m)$ from which the first row is removed is denoted as $(m)_{n-1}$.

The Wigner coefficients  are matrix elements  of an interwinnig operator
$\Phi: V^{[\bar{m}]_{n}}\rightarrow V^{[M]_{n}}\otimes V^{[m]_{n}}$.
All such operators are indexed by tableaux
% $(\Gamma)\in
%V^{[M]_{n}}$,
% такими что $$[\bar{m}]_n=\Delta ((\Gamma))+[m]_n.$$ Этот вектор
%  по оператору   $\Phi$  определяется следующим образом
such that

$$\Phi((m)_{max})=(\Gamma)\otimes (m)_{max}+l.o.t.,\,\,\, (\Gamma)\in V^{[M]_{n}}$$  where $l.o.t.$ (lower order terms)
denotes a sum of tensor products of weight vectors where the second vector has a weight lower than $[m]_n$.

The corresponding Wigner coefficient is denoted as

\begin{equation}\label{viggenn}
<\begin{pmatrix}  [\bar{m}]_n \\  (\bar{m})_{n-1} \end{pmatrix}
\begin{pmatrix}  (\Gamma)_{n-1} \\ [ M]_n\\ (M)_{n-1}\end{pmatrix}
\begin{pmatrix}[m]_n\\ (m)_{n-1}\end{pmatrix}>.
\end{equation}

This coefficient can be  non-zero only if $[\bar{m}]_n=\Delta(\Gamma)+[m]_n$.
The corresponding Wigner coefficient is denoted as
 \begin{align}\label{redviggenn}
<\begin{pmatrix}  [\bar{m}]_n \\  [\bar{m}']_{n} \\
[\bar{m}]_{n-1}
\end{pmatrix} \begin{vmatrix}  (\Gamma)_{n-1} \\ [ M]_n\\ (\gamma)_{n-1}
\end{vmatrix} \begin{pmatrix}[m]_n\\ [m']_{n}  \\  [m]_{n-1}
\end{pmatrix}>
  \end{align}

  \subsection{ Fundamental Wigner coefficients}

% В данном
%разделе мы введён некоторый класс тензорных операторов, или
%коэффициентов Вигнера. Для этого класса отсутствует проблема
%кратности.

In the present paper only the Wigner coefficients are considered for which  $[M]_n=[1,0,...,0]=[1 \dot{0}]_{n}$,  that is when   the tensor factor
$V^{[M]_n}$ is a standard representation. Such Wigner coefficients are called fundamental.

Note that weight vectors $(m)$  of the standard  representations are completely defined by their weights $\Delta(m)=[0,...,\pm 1,...,0]$, where $\pm 1$ occurs at the place $i$. If $i=0$ then only $1$ is allowed. Thus the Wigner coefficient is of type can be denoted as

\begin{equation}
\begin{pmatrix} i \\ [1 \dot{0}]_{n}\\j  \end{pmatrix}
\end{equation}

\section{The method of $Z$-invariants}

\label{oldvector}

Let us  explain the method of $Z$-invariants,  that allows to describe effectively   the space of  $g_{n-1}$-highest vectors in a  $g_n$-representation.

\subsection{ Realization of a representation on the space of functions on uppertriangular matrices.}

Let $g$ be a classical Lie algebra, denote as   $n$ its rank, and let  $G$ be the corresponding Lie group.

Let $$G=Z_{-}DZ_{+}$$be the Gauss decomposition.  Denote the group $Z_{+}$ of upper triangular matrices with units on the diagonal shortly as  $Z$.

Let us construct a realization of a representation with the highest weight $[m_{-n,n},...,m_{-1,n}]$
on the space of polynomial functions on $Z$ .

The elements of the matrix $z\in Z_{}$ are denoted as $z_{ij}$.
Thus the functions are polynomials in variables $z_{ij}$.

Define the function $\alpha$ on the space of diagonal matrices as follows

$$\alpha(\delta)=\delta_{-n}^{m_{-n,n}}...\delta_{-1}^{m_{-1,n}},$$ where
$\delta\in D$  and $\delta_{-n}$,...,$\delta_{-1}$ are diagonal elements of $\delta$.

Define the action $T_g$ of the element $g\in G$ on a function $f(z)$, $z\in
Z$ as follows
\begin{equation}\label{deist}(T_{g}f)(z)=\alpha(\tilde{\delta})f(\tilde{z}),\end{equation} where
$\tilde{\delta}$ and $\tilde{z}$ througn the Gauss decomposition of   $zg$
$$zg=\tilde{\zeta}\tilde{\delta}\tilde{z}, \,\,\,\tilde{\zeta}\in
Z_{-},\,\,\,\, \tilde{\delta}\in D,\,\,\, \tilde{z}\in Z.$$

In the space of function on $Z$ the subspace of functions $f$,
that  form an irreducible representation with the highest weight  $m=[m_{-n,n},...,m_{-1,n}]$
is defined as follows.

Let  $\mathcal{O}_i$ be an operator on the space of functions on $Z$, which is a left infinitesimal  shift on the $i$-root element.  Then the space of functions we are looking for is the solution space of the following system of differential equations which is called the indicator system
\begin{equation}\label{indic}\mathcal{O}_i^{r_i+1}f=0,\,\,\,r_i=\frac{2(m,\omega_i)}{(\omega_i,\omega_i)},\,\,\,i=1,...,n,\end{equation}

where $\omega_i$ are fundamental weights for $g_n$. We call $r_i$ the exponents of the system \eqref{indic}.

\subsection{The highest vectors for the subalgebra  $g_{n-1}\subset
g_n$, that preseves the coordinates $x_{-n},x_{n}$.}\label{ve}

 Let  $G_n=Sp_{2n}$,  $O_{2n+1}$, $O_{2n}$.
Identify $G_{n-1}\subset G_n$ with a subgroup in $G_n$, that preserves coordinates $x_{-1}$, $x_{1}$. The subgroup $Z$ in
$G_n$ is denoted as  $Z_{n}$. Obviously $Z_{n-1}=Z_{n}\cap
G_{n-1}$.  Then $G_{n-1}$-highest vectors correspond to polynomials that are invariant under the action of the subgroup $Z_{n-1}$ and satisfy the indicator system.

Let us find polynomials invariant under the action of  $Z_{n-1}$. In the case $Sp_{2n}$ it is done in \cite{zh}, in the case $O_{2n+1}$,
$O_{2n}$ this can be done in a similar way. The answer is the following.

The polynomials invariant under the action of $Z_{n-1}$ are polynomials in variables $z_{-1,2},...,z_{-1,n}$,
$z_{1,2},...,z_{1,n}$ in the case of orthogonal groups and in variables
$z_{-1,1},z_{-1,2}...,z_{-1,-2}$, $z_{1,2},...,z_{1n}$ in the case of symplectic groups.

%Или же, при альтернативном варианте,  инвариантные под действием
%$Z_{n-1}$ полиномы - это полиномы от переменных
%$z_{-n,-1},...,z_{-2,-1}$, $z_{-n,1},...,z_{-2,1}$ в случае
%ортогональных групп и $z_{-n,-1},...,z_{-2,-1}$,
%$z_{-n,1},...,z_{-2,1}, z_{-1,1}$ в случае симплектической группы.

%%Индикаторная же система в точности совпадает с \eqref{indic}.

\section{The sketch of the construction}
\label{shema}

Let us give a sketch the construction of  the Gelfand-Tsetlin type base,  the calculation of Wigner coefficients  and reduced Wigner coefficients and the derivation of formulas for the action of generators of the algebra.

The induction in  $n$ is used. The subalgebra
$g_{n-1}\subset g_n$ is identified with the subalgebra that preserves the coordinates $x_{-1},x_{1}$.

At first the case $n=1$ is considered, this is the case of algebras
$\mathfrak{o}_3$, $\mathfrak{sp}_2$, $\mathfrak{o}_2$, then the case $n=2$ is considered,
this is the case of algebras $\mathfrak{o}_5$, $\mathfrak{sp}_4$,
$\mathfrak{o}_4$.
% Конструкции
%будут опять различны в этих трёх случаях.

Finally the the passage from $g_{2}$  to  $g_n$ is considered. This passage is done simultaneously for all algebras $B$, $C$, $D$.
In this step the Gelfand-Tsetlin base for $g_n$ is constructed,
Wigner coefficients for  $g_{n-1}$ corresponding to tensor multiplication on the standard representation are obtained.
As before the Wigner coefficients in the case $g_2$  are calculated separately for the three series of algebras and the Wigner coefficients for $g_n$, $n>2$ are calculated in a similar way.
 Finally  the explicit formulas for the action of generator $g_n$ in the base are derived.

Note that the Wigner coefficients for $g_{n-1}$ are calculated when the algebra  $g_n$ is considered.
Let us stress that the passage from $g_{2}$ to $g_{n}$ is very close to the passage from
$\mathfrak{gl}_{3}$ to $\mathfrak{gl}_{n+1}$. This relation plays a key role in the construction.

Is is not necessary to calculate the formulas for the action of all generators. Indeed, it is enough to obtain the formulas for the action of elements  $e_{\pm\alpha}$, where $\alpha
$ are simple roots and also formulas for the action of Cartan elements.
Also let us note that  the matrices that define the action the element $e_{-\alpha}$ can be obtained from the matrix that correspond to  $e_{-\alpha}$ by conjugation,  hence it is enough to  consider only one of the roots $\pm \alpha$.

It is well-know that for the Gelfand-Tsetlin type base the following statement is true. The subalgebra $g_k\subset g_n$ changes only the part of the Gelfans-Tsetlin tableau that corresponds to $g_k$, and this action does not depend on the upper rows.

Thus we obtain that it is enough to calculate the weight of the base vector and also the action of the operators $F_{-1,0}$  for $n=1$,  and $F_{-1,-2}$, $F_{-2,1}$ for  $n=2$
and  $F_{-1,-2}$ for  $n>2$.

%\begin{enumerate}
%\item При рассмотрении случая $n=1$ - оператора $F_{-1,0}$ для
%$\mathfrak{o}_3$.
%\item При рассмотрении случая $n=1$ - операторов $F_{-1,-2}$,
%$F_{-2,1}$.
%\item  При рассмотрении случая $n>2$ - операторов $F_{-1,-2}$.
%\end{enumerate}

\subsection{Notations}

For the rows $[m]_n=[m_{-n,n},...,m_{-1,n}]$ and
$[m]_{n-1}=[m_{-n,n-1},...,m_{-2,n-1}]$ introduce as in
\cite{1968} the following symbols \footnote{In \cite{1968} it is suggested that $m_{-1,n}=0$,
that is why the formula \eqref{redsp2} is slightly different from
\cite{1968}}.

Let

\begin{align}
\begin{split}\label{redsp2}&\begin{vmatrix}n \\
i_1: n-1\end{vmatrix}^{[m]_n,[m]_{n-1}}=\big(
(m_{-2,n-1}-m_{-1,n})\frac{\Pi_{j=-1}^{-n}(m_{j,n}-m_{i_1,n-1}-j+i_1+1)}{\Pi_{j=-2,j\neq
i}^{-n}(m_{j,n-1}-m_{i_1,n-1}-j+i_1+1)}\big)^{\frac{1}{2}},\end{split}\end{align}

be the reduced matrix element for the
 $\mathfrak{gl}_{n-1}$-tensor operator  $E_{n,i}$.
Let

\begin{align}\label{tredsp22}\begin{split}&\begin{vmatrix}i_1: n \\ i_2 : n-1\end{vmatrix}^{[m]_n,[m]_{n-1}}=S(i_2-i_1)\big(
\frac{\Pi_{j=-2,j\neq
i_1}^{-n}(m_{j,n-1}-m_{i_1,n}-j+i_1)}{\Pi_{j=-1,j\neq
i_1}^{-n}(m_{j,n}-m_{i_1,n}-j+i_1)}\cdot \\&\cdot
   \frac{\Pi_{j=-1,j\neq
i_2}^{-n}(m_{j,n}-m_{i_2,n-1}-j+i_2+1)}{\Pi_{j=-2,j\neq
i_2}^{-n}(m_{j,n-1}-m_{i_2,n-1}-j+i_2+1)}\big)^{\frac{1}{2}}
,\end{split},
\end{align}

be  a $\mathfrak{gl}_n$ reduced Wigner coefficient. Here $S(x)=sign(x)$, $S(0)=1$.
Let

\begin{align}
\begin{split}
\label{wigsp222}
&\begin{vmatrix} i: n \\  n-1\end{vmatrix}^{[m]_n,[m]_{n-1}}=\big ( \frac{\Pi_{j=-2}^{-n} (m_{j,n-1}-m_{i,n}-j+i )
 }{ \Pi_{j=-1,j\neq i}^{-n} ( m_{i,n}-m_{i,n}-j+i   )   }
 \big)^{\frac{1}{2}}.
\end{split}
\end{align}

be  $\mathfrak{gl}_{n}$ Wigner coefficient.

%С помощью равенства $m_{-k}=m_{k}$ эти символы определяются также и
%для наборов $[m_{-n},...,m_{-1}]$.

\section{$n=1$: the construction for the algebras $\mathfrak{o}_3$, $\mathfrak{sp}_2$, $\mathfrak{o}_2$.}
\label{step1}

Let us construct the Gelfand-Tsetlin base and derive formulas for the action of generators of these algebras.

For the algebras  $\mathfrak{sp}_2$, $\mathfrak{o}_2$ this is  a trivial task, since these
algebras are one-dimensional and their irreducible representations are also one-dimensional.
They are defined by the highest weight which is just one number $m_{-1}$. Consider the nontrivial case $\mathfrak{o}_3$.

\subsection{The case $\mathfrak{o}_3$}
 The  Gelfand-Tsetlin
base coincides with the  standard weight base in a
  $\mathfrak{sl}_2$-representation.
  But an indexation that we construct differs from the standard
   one.

 The nonstandard indexation is useful in further investigations  since it allows to establish
  a relation for Gelfand-Tsetlin bases for  the algebras of series  $B$  and $A$ for  $n=2$.

Let us use the realization of a representations on the space of
functions on   $Z$.  The indicator system in the case of the highest
weight $m_{-1,1}$ (which is an inter or a half-integer number) is

\begin{equation}
(\frac{\partial }{\partial z_{-1,0}})^{2m_{-1,1}}f=0.
\end{equation}

Thus  the base vectors in the representation in this realization are
monomials $z_{-1,0}^k$, $k=0,...,2m_{-1,1}$. In \cite{zh} it is
shown that the vector corresponding to the monomial  $z_{-1,0}^k$,
 has a weight $m_{-1,1}-k$.

However when the Wigner coefficient for $\mathfrak{sl}_2$ are calculated the orthonomal
base is used, that is the base
$$e_{\mu}=\frac{z_{-1,0}^{m_{-1,1}-\mu}}{\sqrt{(m_{-1,1}-\mu)!(m_{-1,1}+\mu)!}},$$
$\mu=-m_{-1,1},...,m_{-1,1}$ is the weight of the vector.  Below we use this base.

Let us  construct two numbers. Put $m'_{-1,1}=[\frac{k}{2}]$, where
$[ \,.\,]$  is an integer part and let $\sigma_{-1}=0,1$  be the
residue of the division of $k$ by 2. A base vector can be encoded by
a tableau

\begin{align}\begin{split}\label{gco3}
& m_{-1,1}\\
& \sigma_{-1}\,\,m'_{-1,1},\end{split}
\end{align}

The following inequality holds $m_{-1,1}\geq m'_{-1,1}\geq 0$. Also
if $m_{-1,1}=m'_{-1,1}$ and $m_{-1,1}$ is an integer then
$\sigma_{-1}=0$.

The weight of the vector encoded by a tableau can be calculated by
the formula $m_{-1,1}-2m'_{-1,1}-\sigma_1$.

Let us find the action of the operator  $F_{-1,0}$. Put
$g=exp(tF_{-1,0} )\in Z$, then $T_gf(z)=f(zg)$. Under the action of
$F_{-1,0}$ the vector  $z_{-1,0}^k$
   is mapped to $kz_{-1,0}^{k-1}$. One can easily  write how the numbers $m'_{-1,1}$ and $\sigma_{-1}$ change under the action of $F_{-1,0}$.
   This gives us a necessary condition for the matrix element to be non-zero.  In this case the matrix element equals to the reduced
   matrix element of a  $\mathfrak{o}_1$-tensor operator $F_{-1,0}$.
   Thus the following theorem takes place.

Put $[m]_{1}=[m_{-1,1},0]$, $[m']_{1}=[m'_{-1,1}]$.

\begin{thm}

Under the action of $F_{-1,0}$   the tableau $(m)$ changes in the
following way.

    \begin{enumerate}

\item If $\sigma_{-1}=0$, then $\sigma_{-1}$ diminishes by $1$,
$m'_{-1,1}$  remain unchanged,  the resulting tableau is multiplied
to
%$\sqrt{(2m'_{-1,1})(2m_{-1}-2m'_{-1,1}+1)}=
$\begin{vmatrix}2 \\
-2: 1\end{vmatrix}^{[m]_{1},[m']_{1}}$.

\item If  $\sigma_{-1}=1$, then $m'_{-1,1}$ diminishes by $1$,
$\sigma_{-1}$ turns to $0$, the resulting tableau is multiplied to
%$\sqrt{(2m'_{-1,1}+1)(2m_{-1}-2m'_{-1,1})}=
$\begin{vmatrix}2 \\
-2: 1\end{vmatrix}^{[m]_{1},[m']_{1}}$.
\end{enumerate}
\end{thm}

\section{$n=2$:  the construction for the algebras $\mathfrak{o}_5$, $\mathfrak{sp}_4$, $\mathfrak{o}_4$.}
\label{step2}

In the present Section the restrictions $\mathfrak{o}_5\downarrow
\mathfrak{o}_3$, $\mathfrak{sp}_{4}\downarrow\mathfrak{sp}_2$,
$\mathfrak{o}_4\downarrow \mathfrak{o}_2$ are investigated. Then the
Gelfand-Tsetlin base is constructed and  the action of the
generators in this base is investigated.

The construction in the three cases $\mathfrak{o}_5$,
$\mathfrak{sp}_4$, $\mathfrak{o}_4$ are different. For
$\mathfrak{o}_5$ and $\mathfrak{sp}_4$  the corresponding problem of
restriction is solved by establishing a relation with the analogous
problem of restriction  $\mathfrak{gl}_3\downarrow \mathfrak{gl}_1$.

Note that the Wigner coefficients for
$\mathfrak{o}_3=\mathfrak{sl}_2$
 are known thus for $n=2$ it is not necessary to calculate the Wigner coefficients.

\subsubsection{The auxiliary problem of restriction $\mathfrak{gl}_3\downarrow \mathfrak{gl}_1$}

Identify the algebra $\mathfrak{gl}_3$  with the algebra of
endomorphisms of the space with coordinates $x_{-2},x_{-1},x_1$, and
$\mathfrak{gl}_1$ is the subalgebra that preserves the last two
coordinates. The problem of restriction $\mathfrak{gl}_3\downarrow
\mathfrak{gl}_1$ is equivalent to the problem of construction of
weight base in a $\mathfrak{gl}_3$-representation.

We use the realization of a representation with the highest weight
$[\lambda_{-2},\lambda_{-1},\lambda_{1}]$ on the space of polynomials
on the group $Z$, that is on the space of polynomials in variables
$z_{-2,-1},z_{-2,1},z_{-1,1}$. The indicator system is of type

\begin{align}\begin{split}\label{indgl3}
 & (\frac{\partial}{\partial z_{-2,-1}}+z_{-1,1}\frac{\partial}{\partial z_{-2,1}})^{r_{-2}+1}f=0,\\&
(\frac{\partial}{\partial z_{-1,1}})^{r_{-1}+1}f=0, \end{split}
\end{align}

where $r_{-2}=\lambda_{-2}-\lambda_{-1}$, $r_{-1}=\lambda_{-1}$.

The solution space of this system is spanned by polynomials of type
\begin{equation}\label{solgl3}
f=f_0(z_{-2,-1},z_{-2,1},z_{-1,1})z_{-1,1}^{p},\end{equation} where
$f_0$ is a polynomial solution of the first equation, that is not divisible by
$z_{-1,1}$. The following inequality must hold
$$deg_{z_{-1,1}}f_0+p\leq r_{-1},$$

where $deg_{z_{-1,1}}f_0$ is the degree of $f_0$ as a polynomial in $z_{-1,1}$.

%Несложно выписать базис в пространстве полиномиальных решений этой
%системы. Все решения являются линейными комбинациями следующих
%многочленов

%\begin{align}\begin{split}\label{solgl3}
%f_{\mathfrak{gl}_3}=z_{-1,1}^{p_{-1}}(z_{-2,-1}z_{-1,1}-%z_{-2,1})^{p_{-2}}%(z_{-2,-1}+z_{-1,1}z_{-2,1})^{q_{-2}},\\
%p_{-2}+q_{-2}\leq r_{-2},\\
%p_{-1}+p_{-2}+q_{-2}\leq r_{-1}.\end{split}
%\end{align}

There exist a base in the representation whose vectors are encoded by Gelfand-Tsetlin tableaux

\begin{align}\begin{split}\label{gcgl3}
& \lambda_{-2} \,\,\,\, \lambda_{-1}   \,\,\,\, \lambda_{1}\\
& \,\,\,\lambda_{-2,2}\,\,\lambda_{-1,2}\\
&\,\,\,\,\,\,\,\,\lambda_{-2,1}.\end{split}
\end{align}

\subsection{The construction in the case $\mathfrak{o}_5$}
Consider a  representation of $\mathfrak{o}_5$ with the highest weight $[m_{-2,2},m_{-1,2}]$.

\subsubsection{The problem of restriction  $\mathfrak{o}_5 \downarrow \mathfrak{o}_3$}\label{indb}

Using the method of  $Z$-invariants (see Section \ref{ve}) let us investigate how an $\mathfrak{o}_5$-representation branches when we restrict the algebra
$\mathfrak{o}_5 \downarrow \mathfrak{o}_3$.  The polynomials on
$Z_{\mathfrak{o}_5}$, that are  invariant under the action of
$Z_{\mathfrak{o}_3}$, are polynomials in variables $z_{-2,-1},
z_{-2,1}, z_{-1,0}$.

The indicator system is

\begin{align}\begin{split}\label{indo5}&
(\frac{\partial}{\partial z_{-2,-1}}+z_{-1,1}\frac{\partial}{\partial z_{-2,1}})^{r_{-2}+1}f=0,\\&
(\frac{\partial}{\partial z_{-1,0}})^{r_{-1}+1}f=0,\end{split}
\end{align}

where $r_{-2}=m_{-2,2}-m_{-1,2}$, $r_{-1}=2m_{-1,2}$. Note that $z_{-1,1}=-\frac{z_{-1,0}^2}{2}.$

The solution space of this system is spanned by polynomials of type
\begin{equation}\label{solo5}
f=f_0(z_{-2,-1},z_{-2,1},z_{-1,1})z_{-1,0}^{p},\end{equation} where
$f_0$ is a polynomial solution of the first equation that is not divisible by
$z_{-1,1}$. The following inequality must take place
$$2deg_{z_{-1,1}}f_0+p\leq r_{-1}.$$

%Несложно выписать базис в пространстве полиномиальных решений этой
%системы. Все решения являются линейными комбинациями следующих
%многочленов

%\begin{align}\begin{split}\label{solo5}
%f_{\mathfrak{o}_5}=z_{-1,0}^{p_{-1}}(z_{-2,-1}z_{-1,1}-%z_{-2,1})^{p_{-2}}%(z_{-2,-1}+z_{-1,1}z_{-2,1})^{q_{-2}},\\
%p_{-2}+q_{-2}\leq r_{-2},\\
%p_{-1}+2p_{-2}+2q_{-2}\leq r_{-1}.\end{split}
%\end{align}

Let us establish a correspondence between the problems of restriction $\mathfrak{o}_5 \downarrow
\mathfrak{o}_3$ and $\mathfrak{gl}_3\downarrow \mathfrak{gl}_1$.
Consider the cases when $p=2p'$ is even and when   $p=2p'+1$ is odd.

\begin{enumerate}

\item In the case $p=2p'$ to the solution \eqref{solo5} with  the parameter  $p$ of the system \eqref{indo5} with  exponents $r_{-2},r_{-1}$
 there correponds a solution \eqref{solgl3} with the parameter   $p'$ of the system \eqref{indgl3} with exponents $r_{-2},\frac{r_{-1}}{2}$

 $$f_0(z_{-2,-1},z_{-2,1},z_{-1,1})z_{-1,0}^{p}\mapsto f_0(z_{-2,-1},z_{-2,1},z_{-1,1})z_{-1,1}^{[\frac{p}{2}]}. $$

This correspondence is a bijection between the solution space of \eqref{indo5} with an even parameter $p$ and the solution space of
\eqref{indgl3}.

\item In the case $p=2p'+1$ to the solution \eqref{solo5} with  the parameter  $p$ of the system \eqref{indo5} with  exponents $r_{-2},r_{-1}$
 there correponds a solution \eqref{solgl3} with the parameter   $p'$ of the system \eqref{indgl3} with exponents $r_{-2},\frac{r_{-1}}{2}$

 $$f_0(z_{-2,-1},z_{-2,1},z_{-1,1})z_{-1,0}^{p}\mapsto f_0(z_{-2,-1},z_{-2,1},z_{-1,1})z_{-1,1}^{[\frac{p}{2}]}. $$

This correspondence is a bijection between the solution space of \eqref{indo5} with an odd parameter $p$ and the solution space of
\eqref{indgl3}.

\end{enumerate}

\subsubsection{The Gelfand-Tsetlin base for $\mathfrak{o}_5$}\label{gcbaseo5}

In the case $\mathfrak{gl}_{3}$
the indicator system with exponents
$r_{-2},\frac{r_{-1}}{2}$  corresponds to the highest weight $[m_{-2,2},m_{-1,2},0]$.
Base vectors of a representation of $\mathfrak{gl}_{3}$ are encoded by tableaux \eqref{gcgl3}.

Denote as $\sigma_{-1}$ the residue of the division of $p$ by $2$.
Using the correspondence that was established in the previous subsection one obtains that $\mathfrak{o}_3$-highest vectors in a
$\mathfrak{o}_5$-representation are encoded by
$\mathfrak{gl}_{3}$-tableaux   (from which the zero in the upper row is removed) with an additional number $\sigma_{-2}$, that is by tableaux of type

\begin{align}\begin{split}\label{gco5}
& \,\,\,m_{-2,2} \,\,\,\, m_{-1,2}   \,\,\,\, \\
& \sigma_{-2} \,\,\,m'_{-2,2}\,\,m'_{-1,2}\\
&\,\,\,\,\,\,\,\,\,\,\,m_{-2,1}.\end{split}
\end{align}

The following inequalities must take place

\begin{align}
& m_{-2,2}\geq m'_{-2,2}\geq m_{-1,2}\geq m'_{-1,2}\geq 0\\& m'_{-2,2}\geq
m_{-2,1}\geq m'_{-1,2}.
\end{align}

Also if  $m'_{-2,2}=0$, then  $\sigma_{-2}=0$.
Indeed, consider the inequality $2deg_{z_{-1,1}}f_0+p\leq r_{-1},$  and divide it by  two. One obtaines
$deg_{z_{-1,1}}f_0+p'+\sigma_{-2}\leq  \frac{r_{-1}}{2}.$
From here we conclude that if $r_{-1}$ is even (that is when the highest weight
$[m_{-2,2},m_{-1,2}]$ is integer) and $p'$ takes he maximm value then $\sigma_{-2}=0$. The fact that  $p'$ takes the maximum value means that
$m'_{-2,2}=0$ (see explicit formula for the polynomial on $Z$,
corresponding to a $\mathfrak{gl}_3$-tableau in \cite{zh}).

To be able to use this indexation of
$\mathfrak{o}_3$-highest vectors for the construction of the Gelfand-Tsetlin type base we must prove the following Proposition.

\begin{prop}\label{stroka05} The number $m_{-2,1}$ in the tableau \eqref{gco5} is the
$\mathfrak{o}_3$-weight of the vector that is encoded by the tableau.
\end{prop}

The proof is elementary, it can be found in Appendix in Section \ref{sootv}.

Using the  inductive procedure of the construction of the Gelfand-Tsetlin type base we obtain that the vectors of the representation are encoded by  a tableau $(m)$,  that is  obtained by adding to the tableau
\eqref{gco5} the tableau  \eqref{gco3}, that is by  a tableau $(m)$ of type

\begin{align}\begin{split}\label{gco55}
& \,\,\,m_{-2,2} \,\,\,\, m_{-1,2}   \,\,\,\, \\
& \sigma_{-2} \,\,\,m'_{-2,2}\,\,m'_{-1,2}\\
&\,\,\,\,\,\,\,\,\,\,\,m_{-2,1}\\
&\,\,\,\,\,\,\,\,\,\,\,\sigma_{-1} \,\,\,m'_{-2,1}.\end{split}
\end{align}

The elements of this tableau must satisfy the inequalities corresponding to the tableaux
 \eqref{gco5},  \eqref{gco3}.

Let us give formulas for the weight
$[\Delta(m)_{-2},\Delta(m)_{-1}]$ of the vector encoded by  tableau. It is enough to give a formula for
$\Delta(m)_{-1}$.

\begin{prop}
\label{weightofo5}

$\Delta(m)_{-1}=-2\sum_i m'_{i,2}+\sum_i
m_{i,2}+m_{-2,1}-\sigma_{-2}$

\end{prop}

The proof can be found in Appendix in Section
\ref{sootv}.

\subsubsection{Reduced matrix elements for $\mathfrak{o}_5$}
\label{reduct05}

To calculate  matrix elements of generators in the base that we have constructed we need explicit formulas for  reduced matrix elements of operators
 $F_{-1,-2}$ and
$F_{1,-2}$, viewed as $\mathfrak{o}_3$-tensor operators acting between $\mathfrak{o}_3$-representations into which an  $\mathfrak{o}_5$-representation splits

%При вычислении редуцированных матричных элементов можно предполагать,
%что векторы, между %которыми берётся этот элемент, максимальны по отношению к $\mathfrak{o}_3$.

Let $(\bar{m})_{red}$, $(m)_{red}$ be two tableaux of type
\eqref{gco5} that define $\mathfrak{o}_3$-highest vectors in a
$\mathfrak{o}_5$ representation. Let $(\bar{m})^{*}_{red}$ and
$(m)^{*}_{red}$ be two $\mathfrak{gl}_3$-tableux, that are obtained from  $(\bar{m})_{red}$, $(m)_{red}$ by removing the number $\sigma_1$ and adding $0$ to the upper row in the right.

 The reduced matrix element of a $\mathfrak{o}_3$-tensor operator  $F_{\pm 1,-2}$ depends on tableaux  $(\bar{m})_{red}$, $(m)_{red}$. Denote it   $<(\bar{m})_{red}|F_{\pm 1, -2}|(m)_{red}>_{red}$.

Denote as  $<(\bar{m})^{*}_{red}|Е_{\pm 1, -2}|(m)^{*}_{red}>$ a $\mathfrak{gl}_3$-matrix element.

The following lemma takes place.

\begin{lem}\label{gl3o5}  The  matrix elements for $\mathfrak{o}_5$ and
$\mathfrak{gl}_3$ are equal
\begin{equation} <(\bar{m})_{red}|F_{\pm 1, -2}|(m)_{red}>_{red}=<(\bar{m})^{*}_{red}|Е_{\pm 1, -2}|(m)^{*}_{red}>
\end{equation}
\end{lem}

\proof

 The proof uses the following fact. Let
$(m)_{max}$ be a  $\mathfrak{o}_{5}$-tableau that is obtained from
$(m)_{red}$ by adding maximal rows below.

\begin{prop}\label{pp}
The following equality takes place
\begin{align}
<(\bar{m})_{max} \mid F_{\pm 1,-2} \mid (m)_{max}
>=<(\bar{m})_{red} \mid F_{\pm 1,-2} \mid (m)_{red}>_{red}
\end{align}
\end{prop}

The proof can be found in Appendix in Section \ref{redel}.

Let us return to the proof of Lemma. It enough to prove that

$$<(\bar{m})_{max} \mid F_{\pm 1,-2} \mid (m)_{max}
>=<(\bar{m})^{*}_{red}|Е_{\pm 1, -2}|(m)^{*}_{red}>.$$

Let $f=z_{-1,0}^pf_0(z_{-2,-1},z_{-2,1},z_{-1,1})$ be a polynomial corresponding to $(m)_{max}$. One can easily check that under the action of $e^{tF_{-2,-1}}$ it is transformed into the  polynomial
$z_{-1,0}^pf_0(z_{-2,-1}+t,z_{-2,1},z_{-1,1})$.

To the vector $(m)^{*}_{red}$ there corresponds a polynomial
$f^*=z_{-1,1}^{[\frac{p}{2}]}f_0(z_{-2,-1},z_{-2,1},z_{-1,1})$. One can easily check that under the action of $e^{tE_{-2,-1}}$ it is transformed into the polynomial
$z_{-1,1}^{[\frac{p}{2}]}f_0(z_{-2,-1}+t,z_{-2,1},z_{-1,1})$.

Thus the actions  of $F_{-2,-1}$ and $E_{-2,-1}$ on the highest $\mathfrak{o}_3$ and $\mathfrak{gl}_1$  vectors are agreed with  the correspondence that was constructed in Section  \ref{indb}.  Thus the matrix elements of these operators are equal. Hence the matrix elements of the operators $F_{-1,-2}$ and $E_{-1,-2}$ are equal. Using Proposition \ref{pp} we prove the lemma for the operators $F_{-1,-2}$ and $E_{-1,-2}$.

Analogously under the action of the operator $e^{tF_{-2,1}}$, the polynomial  $f$ is transformed into the polynomial
$z_{-1,0}^pf_0(z_{-2,-1},z_{-2,1}+t,z_{-1,1})$. The polynomial  $f^*=z_{-1,1}^{[\frac{p}{2}]}f_0(z_{-2,1},z_{-2,1},z_{-1,1})$ corresponding to $(m)^{*}_{red}$ under the action f the operator $e^{tE_{-2,1}}$ is transformed into the polynomial
$z_{-1,1}^{[\frac{p}{2}]}f_0(z_{-2,-1},z_{-2,1}+t,z_{-1,1})$.

Form this fact we obtain the statement of Lemma for operators $F_{1,-2}$ and $E_{1,-2}$.

\endproof

 Let us find reduced matrix elements for operators $F_{-1,2}$ and $F_{1,2}$ using the lemma.

Put $[m]_{2}=[m_{-2,2},m_{-1,2},0]$, $[m']_{2}=[m'_{-2,2},m'_{-1,2}]$, $[m]_{1}=[m_{-2,1}]$.

Let us prove theorems

\begin{thm}

\label{t1o5}

The reduced matrix element $<(\bar{m})_{red}
\mid F_{1,-2} \mid (m)_{red}>_{red}$ can be nonzero only if the following condition holds. There exist a unique index $i_1=-2$ or $-1$ such that
$\bar{m}'_{i_1,2}=m'_{i_1,2}-1$, and for the other index $i$ one has
$\bar{m}'_{i,2}=m'_{i,2}$. Also $\bar{m}_{-2,1}=m_{-2,1}-1$.

If this condition holds one has
\begin{align}<(\bar{m})_{red} \mid F_{1,-2} \mid (m)_{red}>_{red}=
%\begin{vmatrix}3 \\ i_1,-2 :
%1\end{vmatrix}=
\begin{vmatrix}3 \\ i_1 :
2\end{vmatrix}^{[m]_2,[m']_{2}}\begin{vmatrix}i_1: 2 \\ -2 :
1\end{vmatrix}^{[m']_2,[m]_{1}},
\end{align}

\end{thm}
\proof

By Lemma \ref{gl3o5} it is enough to find  matrix elements for the operator $E_{1,-2}$.
% in the problem of restriction
%$\mathfrak{gl}_{3}\downarrow \mathfrak{gl}_{1}$.

In \cite{1963} in the case $\mathfrak{gl}_{n}$ the following formula is proved

\begin{align}
<(\bar{m}) \mid E_{1,-2} \mid (m)>=\begin{vmatrix}n \\
i_1 : n-1\end{vmatrix}^{[m]_n,[m]_{n-1}}\begin{vmatrix}i_1: n-1 \\ i_2 :
n-2\end{vmatrix}^{[m]_{n-1},[m]_{n-2}}\begin{vmatrix}i_2 : n-2 \\ n-3\end{vmatrix}^{[m]_{n-2},[m]_{n-3}},
\end{align}

where $[m]_k$ is the $k$-th row in the Gelfand-Tsetlin tableau  $(m)$.

It is suggested that one has $\bar{m}_{i_j,j}=m_{i_j,j}-1$, and for $k\neq i_j$
one has $\bar{m}_{i,j}=m_{i,j}$. In all other cases
$<(\bar{m}) \mid E_{1i} \mid (m)>$ vanishes.

If one considered $E_{1,-2}$ as a $\mathfrak{gl}_{n-2}$-tensor operator then the first two factors are a reduced matrix element and the last is a Wigner coefficient.

If one puts $n=3$ then (since  $i_2=-2$)
one gets

\begin{align}
<(\bar{m})_{red} \mid E_{1,-2} \mid (m)_{red}>_{red}=\begin{vmatrix}3 \\
i_1 : 2\end{vmatrix}^{[m]_3,[m]_2}\begin{vmatrix}i_1: 2 \\ -2 :
1\end{vmatrix}^{[m]_2,[m]_1}.
\end{align}

Using now Lemma \ref{gl3o5} we prove the Theorem.

\endproof

%Итак, получена теорема

Let us find the formula for the reduced matrix element of the operator
$F_{-1,-2}$.

\begin{thm}
\label{t2o5}

The reduced matrix element $<(\bar{m})_{red}
\mid F_{-1,-2} \mid (m)_{red}>_{red}$  can be nonzero only if the following condition holds. One has   $\bar{m}_{-2,1}=m_{-2,1}-1$ and all other entries of $(\bar{m})$ and $(m)$ are equal.

 If this condition holds one has

\begin{align}<(\bar{m})_{red} \mid F_{-1,-2} \mid (m)_{red}>_{red}=
\begin{vmatrix}2 \\ -2 : 1\end{vmatrix}^{[m']_{2},[m]_{1}},
\end{align}

\end{thm}

\proof By Lemma \ref{gl3o5} it is enough to find matrix elements $E_{-1,-2}$. One has the formula

\begin{align}
<(\bar{m})_{red} \mid E_{-1-2} \mid (m)_{red}>=\begin{vmatrix}n \\
i_1 : n-1\end{vmatrix}^{[m]_{n},[m]_{n-1}}\begin{vmatrix}i_1: n-1 \\  n-2\end{vmatrix}^{[m]_{n-1},[m]_{n-2}},
\end{align}

where $[m]_k$ is the $k$-th row of the Gelfand-Tsetlin tableau  $(m)$.
It is suggested that $\bar{m}_{i_1,n}=m_{i_1,n}-1$ and
$\bar{m}_{j,n}=m_{j,n}$ for $j\neq i_1$. The first factor is the reduced matrix element. Put $n=2$,  using Lemma \ref{gl3o5} we obtain the statement of the Theorem

\endproof

\subsubsection{Matrix elements}

Using calculation that were done in the previous sections let us derive formulas for matrix elements of generators of $\mathfrak{o}_{5}$ in the base that we have constructed.

It is enough to give formulas for the  matrix elements of operators.
$F_{-1,-2}$ и $F_{1,-2}$.
%Они легко получаются с помощью формул
%Вигнера-Эккарта.
The following theorem takes place.

\begin{thm}

\label{t3o5}

The  matrix element $<(\bar{m})_{red}
\mid F_{1,-2} \mid (m)_{red}>_{red}$ can be nonzero only if the following condition holds. There exist a unique index $i_1=-2$ or $-1$ such that
$\bar{m}'_{i_1,2}=m'_{i_1,2}-1$, and for the other index $i$ one has
$\bar{m}'_{i,2}=m'_{i,2}$. Also $\bar{m}_{-2,1}=m_{-2,1}-1$.

 If this condition holds one has

\begin{align}\begin{split}&
<(\bar{m})| F_{1,-2} |(m)>=
%<(\bar{m})_{red}| F_{1,-2}
%|(m)_{red}>_{red}<(\bar{m})\begin{pmatrix}j\\
%[1\dot{0}]_{n-1} \\ -2 \end{pmatrix} (m)>= \\&=
\begin{vmatrix}3
\\ i_1 : 2\end{vmatrix}^{[m]_2,[m']_2}\begin{vmatrix}i_1: 2 \\ -2 :
1\end{vmatrix}^{[m']_2,[m]_1}\begin{vmatrix}-2 : 2
\\   1  \end{vmatrix}^{[m]_1,[m']_1}
.\end{split}
\end{align}

The matrix element  $<(\bar{m}) \mid
F_{-1,-2} \mid (m)>$ can be nonzero only if the following condition holds. One has
$\bar{m}_{-2,1}=m_{-2,1}-1$ and all other entries of $(\bar{m})$ and $(m)$ are equal.

 If this condition holds one has

\begin{align}\begin{split}&
<(\bar{m})| F_{-1,-2} |(m)>=
%<(\bar{m})_{red}| F_{-1,-2}
%|(m)_{red}>_{red}<(\bar{m})\begin{pmatrix}j\\
%[1\dot{0}]_{n-1} \\ -2 \end{pmatrix} (m)>= \\&=
\begin{vmatrix}2 \\
-2 : 1\end{vmatrix}^{[m']_2,[m]_1}\begin{vmatrix}-2: 2 \\   1  \end{vmatrix}^{[m]_1,[m']_1},
\end{split}
\end{align}

\end{thm}

\proof The Theorem is immediately proved using the Wigner-Ekcart theorem.

\endproof

\subsection{The construction in the case $\mathfrak{sp}_4$}

Consider an irreducible representation of $\mathfrak{sp}_4$  with the highest weight $[m_{-2,2},m_{-1,2}]$. Since the algebra $\mathfrak{sp}_2$
is one dimensional the problem of restriction  $\mathfrak{sp}_4\downarrow
\mathfrak{sp}_2$ is equivalent to the problem of construction of a base in a
$\mathfrak{sp}_4$-representation.

\subsubsection{The problem of restriction  $\mathfrak{sp}_4 \downarrow \mathfrak{sp}_2$}\label{indc}

 Let us investigate the problem of restriction $\mathfrak{sp}_4 \downarrow
\mathfrak{sp}_2$ using the method of $Z$-invariants. The polynomials $Z_{\mathfrak{sp}_4}$ that are invariant under $Z_{\mathfrak{sp}_2}$ are polynomials in variables
$z_{-2,-1}, z_{-2,1}, z_{-1,1}$.

The indicator system is

\begin{align}\begin{split}\label{indsp4}&
(\frac{\partial}{\partial z_{-2,-1}}+z_{-1,1}\frac{\partial}{\partial z_{-2,1}})^{r_{-2}+1}f=0\\&
(\frac{\partial}{\partial z_{-1,1}})^{r_{-1}+1}f=0,\end{split}
\end{align}

where $r_{-2}=m_{-2,2}-m_{-1,2}$, $r_{-1}=m_{-1,2}$.

The system \eqref{indsp4}   is exatly the system \eqref{indgl3}, that appears in the problem of restricton  $\mathfrak{gl}_3\downarrow \mathfrak{gl}_1$.
%, где
%$\mathfrak{gl}_3$-представление имеет  старший вес
%$[m_{-2},m_{-1},0]$.

\subsubsection{The Gelfand-Tsetln type base for $\mathfrak{sp}_4$}\label{gcbasesp4}

In the case $\mathfrak{gl}_{3}$ the indicator system with exponents
$r_{-2},r_{-1}$  corresponds to the highest weight  $[m_{-2,2},m_{-1,2},0]$.
Base vectors of a $\mathfrak{gl}_{3}$-representation are encoded by tableaux  \eqref{gcgl3}. Hence the base vectors of a
$\mathfrak{sp}_4$-representation can be encoded by a $\mathfrak{gl}_3$-tableau (form which we remove zero in the upper row), that is by tableaux of type

\begin{align}\begin{split}\label{gcsp4}
& \,\,\,m_{-2,2} \,\,\,\, m_{-1,2}   \,\,\,\, \\
&\,\,\,\, \,\,\,m'_{-2,2}\,\,m'_{-1,2}\\
&\,\,\,\,\,\,\,\,\,\,\,m_{-2,1}.\end{split}
\end{align}

The entries of this tableau must satisfy the inequalities

\begin{align}
& m_{-2,2}\geq m'_{-2,2}\geq m_{-1,2}\geq m'_{-1,2}\geq 0\\& m'_{-2,2}\geq
m_{-2,1}\geq m'_{-1,2}.
\end{align}

The following Proposition takes place

\begin{prop}\label{strokasp4} The number $m_{-2,1}$ in the tableau \eqref{gcsp4} is the
$\mathfrak{sp}_4$-weight of the  corresponding vector
\end{prop}

The proof is analogous to the proof of the Proposition
 \ref{stroka05}.

Let us  give the formula for the weight
$[\Delta(m)_{-2},\Delta(m)_{-1}]$ of the vector encoded by the tableau $(m)$. It is enough to consider
$\Delta(m)_{-1}$.

\begin{prop}
\label{weightofsp4}

$\Delta(m)_{-1}=-2\sum_i m'_{i,2}+\sum_i m_{i,2}+m_{-2,1}$

\end{prop}

The proof of this Proposition is analogous to the proof of Proposition \ref{weightofsp4}.

\subsubsection{Matrix elements for $\mathfrak{sp}_4$}
 Matrix elements for
$F_{-1,-2}$ and $F_{1,-2}$ can be calculeted if we consider these opearators as
$\mathfrak{sp}_2$-tensor operators acting between
$\mathfrak{sp}_2$-representations into which a
$\mathfrak{sp}_4$-representation splits.

Let us we write the Wigner-Eckart theorem for $F_{-1,-2}$
and $F_{1,-2}$. Since the algebra $\mathfrak{sp}_2$is one dimensional, the  corresponding Wigner coefficient is zero. Hence the matrix elements of $F_{-1,-2}$ and $F_{1,-2}$ equal to the reduced matrix elements.

%Рассуждения полностью аналогичны рассуждениям из раздела
%\ref{reduct05} для случая $\mathfrak{o}_5$.

As in Section \ref{reduct05} let us find the relation of these reduced matrix elements with the analogous matrix elements for $\mathfrak{gl}_3$.

Let $(\bar{m})$, $(m)$  be two tableaux of type \eqref{gcsp4}, and let
$(\bar{m})^{*}$, $(m)^{*}$ be two tableaux for
$\mathfrak{gl}_3$,
  that are obtained from  $(\bar{m})$, $(m)$  by adding zero $0$ to the upper row in the right.
The following Lemma takes place.

\begin{lem}\label{gl3sp4} Matrix elements  for $\mathfrak{sp}_4$ and
$\mathfrak{gl}_3$ are equal
\begin{equation} <(\bar{m})|F_{\pm 1, 2}|(m)>=<(\bar{m})^*|Е_{\pm 1, 2}|(m)^*>
\end{equation}
\end{lem}

The proof is analogous to the proof of Lemma \ref{gl3o5}.

Introduce notations $[m]_2=[m_{-2,2},m_{-1,2},0]$, $[m']_2=[m'_{-2,2},m'_{-1,2}]$, $[m]_{-1}=[m_{-2,1}]$.

Using Lemma \ref{gl3sp4},let us calculate matrix elements for $F_{-1,2}$ and
 $F_{1,2}$. They are given by the follolwing theorem.

\begin{thm}
The matrix element $<(\bar{m}) \mid F_{1,-2} \mid
(m)>$ can be nonzero only if the following condition holds. For one index $i_1=-2$
or $-1$ one has $\bar{m}'_{i_1,2}=m'_{i_1,2}-1$, and for another one has
$\bar{m}_{ij}=m_{ij}$. Also $\bar{m}_{2,1}=m_{2,1}-1$. If this condition holds one has

\label{t1sp4}
\begin{align}<(\bar{m})\mid F_{1,-2} \mid (m)>=
%\begin{vmatrix}3 \\ i_1,-2 :
%1\end{vmatrix}=
\begin{vmatrix}3 \\ i_1 :
2\end{vmatrix}^{[m]_2,[m']_2}\begin{vmatrix}i_1: -2 \\ -2 : 1\end{vmatrix}^{[m']_2,[m]_1},
\end{align}

The matrix element $<(\bar{m}) \mid F_{-1,-2} \mid
(m)>$ can be nonzero only if the following condition holds. One has
$\bar{m}_{-2,1}=m_{-2,1}-1$ and all othe entries of $(\bar{m})$ and $(m)$ are equal. If this condition holds one has

\begin{align}<(\bar{m}) \mid F_{-1,-2} \mid (m)>=
\begin{vmatrix}2 \\ -2 : 1\end{vmatrix}^{[m']_2,[m]_1},
\end{align}

\end{thm}

The proof of this Theorem is base on Lemma \ref{gl3sp4} and is analogous to the proof of Theorem \ref{t1o5}.

\subsection{The construction in the case $\mathfrak{o}_4$}

Let us be given a  $\mathfrak{o}_4$-represention with the highest weight  $[m_{-2,2},m_{-1,2}]$.  Since the algebra $\mathfrak{o}_2$
is one dimensional the problem of restriction $\mathfrak{o}_4\downarrow
\mathfrak{o}_2$ is equivalent to the problem of construction of a base in a
$\mathfrak{o}_4$-representation.

\subsubsection{The problem of restriction $\mathfrak{o}_4 \downarrow \mathfrak{o}_2$}\label{indd}

Consider the problem of restriction
$\mathfrak{o}_4 \downarrow \mathfrak{o}_2$. The polynomials on
$Z_{\mathfrak{o}_4}$ that are invariant under
$Z_{\mathfrak{o}_2}$ are the polynomials in variables $z_{-2,-1},
z_{-2,1}$.

The indicator system is

\begin{align}\begin{split}\label{indo4}
(\frac{\partial}{\partial z_{-2,-1}})^{r_{-2}+1}f=0,\\
(\frac{\partial}{\partial z_{-2,1}})^{r_{-1}+1}f=0,\end{split}
\end{align}

where $r_{-2}=m_{-2,2}-m_{-1,2}$, $r_{-1}=m_{-2,2}+m_{-1,2}$.

All polynomial solutions are linear combination of the following polynomials (we use the standard normalization for
$\mathfrak{sl}_2$)

\begin{align}\begin{split}\label{solo4}
& f=\frac{z_{-2,-1}^{p}}{\sqrt{p!(2r_{-1}-p)!}}\frac{ z_{-2,1}^{q}}{\sqrt{q!(2r_{-2}-q)!}},\\
& 0\leq p\leq r_{-2},\\
& 0\leq q \leq r_{-1}.\end{split}
\end{align}

\subsubsection{The Gelfand-Tsetlin type base for $\mathfrak{o}_4$}\label{gcbaseo4}

 We see that there exist a base whose vectors are defined by two numbers $p$, $q$. Following Molev (see \cite{M})
let us define another indexation. Put

\begin{equation}
m'_{-2,2}=m_{-2,2}-min\{p,q\},\,\,\,
m_{-2,1}=m_{-1,2}-p+q.
\end{equation}

\begin{prop}\label{pq}
The numbers $p$ and $q$ are reconstructed as follows

\begin{enumerate}
\item Если $m'_{-2,1}-m_{-1,2}>0$, то $p=m_{-2,2}-m'_{-2,2}$, $q=m_{-2,2}-m'_{-2,2}+m'_{-2,1}-m_{-1,2}$.
\item Если $m'_{-2,1}-m_{-1,2}\leq 0$, то $q=m_{-2,2}-m'_{-2,2}$, $p=m_{-2,2}-m'_{-2,2}-m'_{-2,1}+m_{-1,2}$.
\end{enumerate}
\end{prop}

Let us construct the Gelfand-Tsetlin type tableau in the following way

\begin{align}\begin{split}\label{gco4}
& \,\,\,m_{-2,2} \,\,\,\, m_{-1,2}   \,\,\,\, \\
&\,\,\,\,\,\,  \,\,\,m'_{-2,2}\\
&\,\,\,\,\,\,\,\,\,\,\,m_{-2,1}\end{split}
\end{align}

The equalities \eqref{solo4}  are satisfied if and only if (see \cite{M}), when the following inequalities take place

\begin{align}
& m_{-2,2}\geq m'_{-2,2}\geq |m_{-1,2}|\\
& m'_{-2,2}\geq |m_{-2,1}|
\end{align}

One has

\begin{prop}\label{stroka04} The number $m_{-2,1}$ in the tableau \eqref{gco4} is the
$\mathfrak{o}_2$-weight of the corresponding vector.
\end{prop}

The proof is analogous to the proof of Proposition \ref{stroka05}.

Let us give formulas for the weight
$[\Delta(m)_{-2},\Delta(m)_{-1}]$ of the vector encoded by the tableau $(m)$. It is enough to consider the component
$\Delta(m)_{-1}$.

\begin{prop}
\label{weightofo4}

$\Delta(m)_{-1}=-2\sum_i m'_{i,2}+\sum_i m_{i,2}+m_{-2,1}$

\end{prop}

The proof is analogous to the proof of Proposition
\ref{weightofo5}.

\subsubsection{Matrix elements}

Since the construction of the Gelfand-Tsetlin type base for
$\mathfrak{o}_4$ has another character than the construction for the algebras
$\mathfrak{o}_5$ and $\mathfrak{sp}_4$ the calculation of matrix elements does not use the Wigner-Eckart theorem. It is based on the isomorphism $\mathfrak{o}_4\simeq
\mathfrak{sl}_2\oplus\mathfrak{sl}_2$.

The operator $F_{-1,-2}$ increases the number $p$  by $1$  and multiplies the vector on the reduced matrix element
%$\begin{vmatrix} -2\\-2:-1\end{vmatrix}$.
$\sqrt{p_{}!(2r_{-2}-p_{})!}$.

The operator $F_{-1,2}$ increases the number  $q$ by $1$ and multiplies the vector on the reduced matrix element
%$\begin{vmatrix} -2\\-2:-1\end{vmatrix}$.
$\sqrt{q_{}!(2r_{-1}-q_{})!}$.

When we pass to the Gelfand-Tsetlin type tableaux then using Proposition \ref{pq}, we obtain the theorem

\begin{thm}

The action of $F_{-1,-2}$ on the tableau \eqref{gco4} is described as follows. If $m'_{2,1}-m_{-1,2}+1\leq 0$, then $F_{-1,-2}$ diminishes $m_{-2,1}$ by  $1$. If  $m'_{2,1}-m_{-1,2}+1> 0$ then also $m'_{-2,2}$ diminishes by $1$.
In both cases the vector is multiplied onto
$\sqrt{p_{}!(2r_{-2}-p_{})!}$ (see Proposition \ref{pq}).

The action of $F_{-1,2}$ on the tableau \eqref{gco4} is described as follows. If  $m'_{2,1}-m_{-1,2}-1\geq 0$, then $F_{-1,-2}$ increases $m_{-2,1}$ by $1$. If $m'_{2,1}-m_{-1,2}+1< 0$  then also  $m'_{-2,2}$ diminishes by $1$.
In both cases the vector is multiplied by
$\sqrt{q_{}!(2r_{-1}-q_{})!}$ (see Proposition  \ref{pq}).

\end{thm}

% Используя вычисления, проведённые в трех
%предыдущих подразделах, выведем формулы для матричных элементов
%действия алгебры $\mathfrak{o}_{4}$ в построенном базисе.

%Достаточно привести формулы  действия $F_{-1,-2}$ и $F_{1,-2}$. С
%помощью теоремы Вигнера-Эккарта доказывается следующая теорема.

%\begin{thm}

%\label{t3o4}
%\begin{align}\begin{split}&
%<(\bar{m})| F_{-1,-2} |(m)>=
%%<(\bar{m})_{red}| F_{-1,-2}
%%|(m)_{red}>_{red}<(\bar{m})\begin{pmatrix}j\\
%%[1\dot{0}]_{n-1} \\ -2 \end{pmatrix} (m)>= \\&=
%\begin{vmatrix}2 \\
%-2 : 1\end{vmatrix}\begin{vmatrix}-2: 2 \\   1  \end{vmatrix},
%\end{split}
%\end{align}

%%где предполагается, что $\bar{m}_{i_1,n}=m_{i_1,n}-1$ и
%%$\bar{m}_{j,n}=m_{j,n}$ при $j\neq i_1$. Первый сомножитель и есть
%%искомый редуцированнй матричный элемент.
%%\end{thm}

%%Аналогично,имеем

%%\begin{thm}
%\begin{align}\begin{split}&
%<(\bar{m})| F_{1,-2} |(m)>=
%%<(\bar{m})_{red}| F_{1,-2}
%%|(m)_{red}>_{red}<(\bar{m})\begin{pmatrix}j\\
%%[1\dot{0}]_{n-1} \\ -2 \end{pmatrix} (m)>= \\&=
%\begin{vmatrix}3
%\\ i_1 : 2\end{vmatrix}\begin{vmatrix}i_1: 2 \\ -2 :
%1\end{vmatrix}\begin{vmatrix}i_2: 2
%\\   1  \end{vmatrix}
%,\end{split}
%\end{align}

%где предполагается, что $\bar{m}_{i_1,1}=m_{i_1,1}-1$ и
%$\bar{m}_{i,n-1}=m_{i,n-1}$ при $i\neq i_1$.
%Первый сомножитель и есть
%искомый редуцированнй матричный элемент.
%\end{thm}

\section{The construction in the case $g_n$}
\label{step3} In the present section in the case $n>2$ we construct a Gelfand-Tsetlin type base in a $g_{n}$-representation, the construction is similar for all algebras. Then we calculate the Wigner coefficients for $g_{n-1}$, the calculations for $g_2$ are different for series $B,C,D$ but the calculations for $g_n$, $n>2$ are similar for all algebras. Then we derive  explicit formulas for the action of generators in the Gelfand-Tsetlin type base.

 Consider a  $g_n$-representation with the highest weight
$[m_{-n,n},...,m_{-1,n}]$.

\subsection{The indicator system in the problem of restriction $g_n\downarrow g_{n-1}$}\label{indn}

The indicator system $I_n$ can be presented as a union of two systems of equation: $I'_n$ and $I_2$. Нere  $I_2$ is an indicator system that appears in the problem of restriction $g_2\downarrow g_1$, it is different for series  $B$, $C$, $D$  (see Sections \ref{indb},
\ref{indc}, \ref{indd}), and the system $I'_n$ is the same for all algebras, this is the following system

\begin{align}\begin{split}
&(z_{-n+1,-1}\frac{\partial }{\partial z_{-n,-1}}+z_{-n+1,1}\frac{\partial }{\partial z_{-n,1}})^{r_{-n}+1}f=0,\\
&...,\\
&(z_{-2,-1}\frac{\partial }{\partial z_{-3,-1}}+z_{-2,1}\frac{\partial }{\partial z_{-3,1}})^{r_{-3}+1}f=0.\end{split}
\end{align}

%где подсистема $I_2$ есть индикаторная система, возникающая в задаче ограничения %$g_{2}\downarrow g_{1}$,  которая была рассмотрена в разделах \ref{indb},
% \ref{indc}, \ref{indd}. Подсистему, образованную уравнениями, не
% входящими в $I_2$ обозначим $I'_n$, таким образом, $I_n=I'_n\cup
% I_2$.

The exponents $r_{i}$ are defined by formulas
\begin{equation}
r_{-n}=m_{-n,n}-m_{-n+1,n},...,r_{-3}=m_{-3,n}-m_{-2,n}.
\end{equation}

The indicator system has this type  not only in the problem of restriction $g_n\downarrow g_{n-1}$
for all algebras  $g_n=\mathfrak{o}_{2n+1},\mathfrak{sp}_{2n},\mathfrak{o}_{2n}$,
but also in the problem of restriction $\mathfrak{gl}_{n+1}\downarrow \mathfrak{gl}_{n-1}$.

Using this remark let us construct a base in the solution space of $I_{n}$, starting from a base in the solution space of
$I_2$. Firstly let us show how this base is constructed in the case
$\mathfrak{gl}_{n+1}$,
and then let us show that this construction is valid also in the case $g_n$.

Consider the algebra  $\mathfrak{gl}_{n+1}$. Let this algebra act in the space with coordinates
 $x_{-n},x_{-n+1},...,x_{-1},x_{1}$. Then $\mathfrak{gl}_{n-1}$-highest vector are indexed by tableaux

\begin{align}\begin{split}\label{dcgln1}
&m_{-n,n},\,\,\,\,m_{-n+1,n},\,\,\, m_{-n+2,n},,....\,\,\,\,\,\,\,\,\,m_{-3,n},\,\,\,\, m_{-2,n},\,\,\,\, m_{-1,n},\,\,\,\, m_{1,n}\\
&\,\,\,\,m_{-n,n-1},\,\,\, m_{-n+1,n-1},,....\,\,\,\,\,\,\,\,\,\,\,\,\, m_{-3,n-1},\,\,\,\, m_{-2,n-1},\, m_{-1,n-1},\,\,\,\, \\
&\,\,\,\,\,\,\,\,m_{-n,n-2},\,\,\, m_{-n+1,n-2},,....\,\,\,\,\,\,\,\,\,\,\,\,\, m_{-3,n-2},\,\,\,\, \,\,\,\,\,m_{-2,n-2},\,\,\,\,\,\,\,\,
\end{split}\end{align}
the entries of this tableau must satisfy the betweeness conditions.

%Let us describe all solutions of $I'_n$ in the polynomial ring
%$\mathfrak{k}[z_{-3,-1},z_{-3,1},...,z_{-n,-1},z_{-n,1}]$, where the field
%$\mathfrak{k}$ is $\mathbb{C}(z_{-2,-1},z_{-2,1},z_{-1,1})$.

The procedure of construction of solutions consists of two steps.

{\bf Step 1.} Let us start with a solution $f(z_{-2,-1},z_{-2,1},z_{-1,1})$ of the system $I_2$,
corresponding to a  $\mathfrak{gl}_3$-tableau

\begin{align}\begin{split}\label{gl3f}
& m_{-2,n},\, m_{-1,n},\, m_{1,n}\\
& \,\,\,\,\,\,\,m_{-2,n-1},\, m_{-1,n-1}\,, \\
& \,\,\,\,\,\, \,\,\,\,\,\, m_{-2,n-2},
\end{split}
\end{align}

Also $f$ is a solution of the sysytem $I_n$, which is encoded by a tableau

\begin{align}
\begin{split}
&m_{-n,n},\,\,\,\,m_{-n+1,n},\,\,\, m_{-n+2,n},,....\,\,\,\,\,\,\,\,\,m_{-3,n},\,\,\,\, m_{-2,n},\,\,\,\, m_{-1,n},\,\,\,\, m_{1,n}\\
&\,\,\,\,m_{-n,n},\,\,\, m_{-n+1,n},,....\,\,\,\,\,\,\,\,\,\,\,\,\, m_{-3,n},\,\,\,\, m_{-2,n-1},\, m_{-1,n-1},\,\,\,\, \\
&\,\,\,\,\,\, \,\,\,\,\,\,\,\,m_{-n,n},\,\,\, m_{-n+1,n},,....\,\,\,\,\,\,\,\,\,\,\,\,\, m_{-3,n},\,\,\,\, \,\,\,\,\,m_{-2,n-2},\,\,\,\,\,\,\,\,
\end{split}
\end{align}

{\bf Step 2.} Let us transform the  solution $f$ of the system $I_n$  to the solution that corresponds to the tableau \eqref{dcgln1}. In \cite{zh} it is shown that the polynomial corresponding to the tableau   \eqref{dcgln1} is of type 

\begin{equation}
\prod_{i=-n}^{-2}
\nabla_{-1,i}^{m_{i,n-1}-m_{i,n-2}}\prod_{i=-n}^{-1} z_{i,1}^{m_{i,n}-m_{i,n-1}}
\end{equation}

where

\begin{align}
\begin{split}\label{om2}
\nabla_{k,i}=\sum_{j_1<...< j_s} c_{j_1...j_s}
E_{k,j_1}E_{j_1,j_2}...E_{j_s,j},\\
c_{j_1,....,j_s}=\prod_{k>j>i,  j\neq j_k}(E_{i,i}-E_{j,j}+j-i).
\end{split}
\end{align}

One can easily check that  the operators $\nabla_{-1,i}$ and the operator of multiplication on   $z_{j,1}$ commute for $i<j$, thus the polynomial corresponding to  \eqref{dcgln1} can be written as

\begin{equation}\label{om3}
(\prod_{i=-n}^{-3}
\nabla_{-1,i}^{m_{i,n-1}-m_{i,n-2}}
  \prod_{i=-n}^{-3}z_{i,1}^{m_{i,n}-m_{i,n-1}})(\nabla_{-1,-2}^{m_{-2,n-1}-m_{-2,n-2}}z_{-2,1}^{m_{-2,n}-m_{-2,n-1}}z_{-1,1}^{m_{-1,n}-m_{-1,n-1}}).
\end{equation}

The second expression in  \eqref{om3} is the polynomial $f$ up to multiplication on a constant. Introduce a notation for the first factor.

\begin{equation}\label{om1}
\Omega^{\mathfrak{gl}_{n+1}}=\prod_{i=-n}^{-3}
\nabla_{-1,i}^{m_{i,n-1}-m_{i,n-2}}
  \prod_{i=-n}^{-3}z_{i,1}^{m_{i,n}-m_{i,n-1}}.
\end{equation}

Thus to the tableau  \eqref{dcgln1} there corresponds the vector

\begin{equation}\label{f52}
\Omega^{\mathfrak{gl}_{n+1}}f.
\end{equation}

Let us fomulate the statement that we have proved. For this let us give a definition.

\begin{defn}
Let  $f$ be a polynomial corresponding to the  tableau \eqref{gl3f}. A set of operators of type \eqref{om1} is called admissible, if after application of operators from this set to  $f$ one obtaines all solutions of
$I'_n$ in the space of polinomials 
$F(z_{-n,-1},z_{-n,1},....,z_{-3,1},z_{-2,-1},z_{-2,1})$ with coeffisients in
$\mathbb{C}$ with the initial condition
$F(0,0,....,0,z_{-2,-1},z_{-2,1})=f(z_{-2,-1},z_{-2,1})$. And also $\Omega_1 f\neq \Omega_2 f$, where
 $\Omega_1\neq \Omega_2$ are admissible operators.
\end{defn}

The previous discussion shows that the set of admissible  operators depend on the elements $m_{-2,n},m_{-2,n-1},m_{-2,n-2}$ of the tableau \eqref{gl3f} only.  An operator $\Omega^{\mathfrak{gl}_{n+1}}$ is admissible if and only if the nubmers $m_{i,n},m_{i,n-1},m_{i,n-2}$ satisfy 

 \begin{align}\begin{split}\label{ineq1}
& m_{-n,n}\geq m_{-n,n-1}\geq  m_{-n+1,n}\geq....\geq m_{-3,n}\geq m_{-3,n-1}\geq m_{-2,n}\\
& m_{-n,n-1}\geq m_{-n,n-2}\geq  m_{-n+1,n-1}\geq .... \geq m_{-3,n-1} \geq m_{-3,n-2} \geq m_{-2,n-1},
% m_{-n,n-2}\geq  m_{-n+1,n-2}\geq ,.... \geq m_{-3,n-2} \geq
% m_{-2,n-2}
 \end{split}
 \end{align}
(
Using explicit formula for the polynomial  $f$ (see   \cite{zh}) corresponding to the tableau   \eqref{gl3f},   one can see that the numbers $m_{-2,n},m_{-2,n-1},m_{-2,n-2}$  are defined by the degrees of variables  $z_{-2,-1},z_{-2,1}$ in the polynomial  $f$.

%  что множество допустимых операторов
%$\Omega^{\mathfrak{gl}_{n+1}}$ зависит от чисел
%$m_{-2,n},m_{-2,n-1}, m_{-2,n-2}$, входящих в диаграмму $D$.  На
%языке многочленов это означает, что множество допустимых операторов
%зависит от степени $f$, как многочлена от $z_{-2,-1},z_{-2,1}$.

Put $z_{-1,1}=1$, one gets the proposition.

\begin{prop}\label{p3r}
In the case $\mathfrak{gl}_{n+1}$ all solution of the system
$I'_n$ in the space of polynomias of type
$F(z_{-n,-1},z_{-n,1},....,z_{-3,1},z_{-2,-1},z_{-2,1})$ with coefficients in
$\mathbb{C}$ with the initial condition
$F(0,0,....,0,z_{-2,-1},z_{-2,1})=f(z_{-2,-1},z_{-2,1})$.
%кодируются таблицами
%\begin{align}\begin{split}\label{tristroki}
%&m_{-n,n},\,\,\,\,m_{-n+1,n},\,\,\, m_{-n+2,n},,....\,\,\,\,\,\,\,\,\,m_{-3,n},\,\,\,\, \\
%&\,\,\,\,m_{-n,n-1},\,\,\, m_{-n+1,n-1},,....\,\,\,\,\,\,\,\,\,\,\,\,\, m_{-3,n-1},\,\,\,\, \, \,\,\,\, \\
%&\,\,\,\,\,\,\,\,m_{-n,n-2},\,\,\,
%m_{-n+1,n-2},,....\,\,\,\,\,\,\,\,\,\,\,\,\, m_{-3,n-2},\,\,\,\,
%\,\,\,\,\,\,\,\,\,\,\,\,\,
%\end{split}\end{align}
%
%элементы которых удовлетворяют условиям промежуточности. В явном
%виде эти решения
can be written as $\Omega^{\mathfrak{gl}_{n+1}}f$, where the set 
of admissible operators $\Omega^{\mathfrak{gl}_{n+1}}$ is defined by the monomials of $f$ or by inequalities 
\eqref{ineq1}.

\end{prop}

Let us give an analogous construction in the case $g_n$.

Put $\mathfrak{k}=\mathbb{C}(z_{-1,0})$   in the case
$\mathfrak{o}_5$, put $\mathfrak{k}=\mathbb{C}(z_{-1,1})$   in the case
$\mathfrak{sp}_4$,  put $\mathfrak{k}=\mathbb{C}$  in the case
$\mathfrak{o}_4$.

 Note that in the case $g_{n}$ the operator $F_{ij}$ acts on a polynomial in variables $z_{-2,-1},z_{-2,1},....,z_{-n,-1},z_{-n,1}$ exactly in the same way as the operator $E_{ij}$ acts on this polynomial.

\begin{defn}
Define a operator $\Omega^{g_{n}}$ by formulas \eqref{om1}, \eqref{om2}, where $E_{ij}$ is replaced to
  $F_{ij}$.
\end{defn}

Now let us construct the set of admissible operators  $\Omega^{g_{n}}$. 
 Let us be given a $g_2$-tableau  $D$ and denote the corresponding polynomial as  $f$.  
Our next pupose is to describe the set of admissible operators  $\Omega^{g_{n}}$, such that applying them to
  $f$ we obtain all solutions of $I'_n$ with the initial consdition $f$ as in Proposition \ref{p3r}.

In the cases $\mathfrak{o}_{2n+1}$, $\mathfrak{sp}_{2n}$ in Sections \ref{indb},  \ref{indc} a correspondence between problems of restriction $\mathfrak{gl}_3\downarrow\mathfrak{gl}_1$  and $\mathfrak{o}_5\downarrow\mathfrak{o}_3$, $\mathfrak{sp}_4\downarrow\mathfrak{sp}_2$
 was established.  This correpondence is generated by a corresponce between solution spaces of indicatir systems. The later correspondence preserves the degrees of variables  $z_{-2,-1},z_{-2,1}$. 
 
The system $I'_n$ is the same in the case $g_n$ and in the case $\mathfrak{gl}_{n+1}$.  Thus one gets that the condition providing that an operator is admissible  in the cases $\Omega^{\mathfrak{o}_{2n+1}}$, $\Omega^{\mathfrak{sp}_{2n}}$   and in the case  $\Omega^{\mathfrak{gl}_{n+1}}$ are the same.  This condition is  the set of inequalities   \eqref{ineq1}.

Now let us find the admissible operators in the case $\mathfrak{o}_{2n}$. In the case  $\mathfrak{gl}_3$ the element $m_{-2,n-1}$ of the Gelfand-Tsetlin tableau can be charateried as follows. We act on the vector corresponding to the tableau by the raising operator $E_{-1,-2}^k$ for the maximum possible  $k$. We obtain a  $\mathfrak{gl}_2$-highest vector. Then $m_{-2,n-1}$ is a $(-2)$-component of it's highest weight.

Consider a vector $\mathfrak{gl}_3$-representation defind by a polynomial  $z_{-2,-1}^pz_{-2,1}^q$, it is a linear combination of tableaux. Let us give a formula for the biggest element $m_{-2,n-1}$ of these tableaux. Under the action
 $E_{-1,-2}^k$ this polynomial is tranformed into  $constz_{-2,-1}^{p-k}z_{-2,1}^{q-k}$. Thus the maximum value of  $k$ is $min(p,q)$. Hence
 \begin{equation} \label{mpro4}  m_{-2,n-1}=m_{-2,n}-min(p,q).\end{equation}

In the case $\mathfrak{o}_4$ to a Gelfand-Tsetlin type tableax there corresponds a polynomial
$const z_{-2,-1}^pz_{-2,1}^q$. The formula for the element $m'_{-2,2}$ of this tableau is just \eqref{mpro4}.

 Also in the case of the problem of restriction $\mathfrak{gl}_3\downarrow \mathfrak{gl}_1$ the element $m_{-2,n-2}$ is the weight of the corresponding vector relatively the subalgebra $\mathfrak{gl}_1$. In the case $g_{n}=\mathfrak{o}_4$ and the problem of the restriction $g_n\downarrow g_{n-1}$ the element $m_{-2,1}$ is the weight relatively $g_{n-1}=\mathfrak{o}_2$.

Thus   in the case $\Omega^{\mathfrak{o}_{2n}}$  and in the case  $\Omega^{\mathfrak{gl}_{n+1}}$ are the same.  This condition is  the set of inequalities   \eqref{ineq1}.

The following analog of the Proposition \ref{p3r} takes place.
% для случая, когда
%рассматриваются полиномы с коэффициентами в поле  $\mathfrak{k}$.

\begin{prop}\label{p33r}
In the case $g_{n}$  all solutions of the system $I'_n$ in the space of polynomials $F(z_{-n,-1},z_{-n,1},....,z_{-3,1},z_{-2,-1},z_{-2,1})$
with coefficients in $\mathfrak{k}$ with the initial condition
$F(0,0,....,0,z_{-2,-1},z_{-2,1})=f(z_{-2,-1},z_{-2,1})$.
%кодируются таблицами
%\begin{align}\begin{split}\label{tristroki}
%&m_{-n,n},\,\,\,\,m_{-n+1,n},\,\,\, m_{-n+2,n},,....\,\,\,\,\,\,\,\,\,m_{-3,n},\,\,\,\, \\
%&\,\,\,\,m_{-n,n-1},\,\,\, m_{-n+1,n-1},,....\,\,\,\,\,\,\,\,\,\,\,\,\, m_{-3,n-1},\,\,\,\, \, \,\,\,\, \\
%&\,\,\,\,\,\,\,\,m_{-n,n-2},\,\,\,
%m_{-n+1,n-2},,....\,\,\,\,\,\,\,\,\,\,\,\,\, m_{-3,n-2},\,\,\,\,
%\,\,\,\,\,\,\,\,\,\,\,\,\,
%\end{split}\end{align}
%
%элементы которых удовлетворяют условиям промежуточности. В явном
%виде эти решения
can be writtem as $\Omega^{g_{n}}f$, where only the admissible
operators $\Omega^{g_{n}}$  are taken.

\end{prop}

Note that if $f$ is a polynomial over $\mathbb{C}$ then
$\Omega^{g_{n}}f$ is a polynomial over $\mathbb{C}$.

%во всех случаях $\mathfrak{gl}_{3}$, $\mathfrak{o}_{5}$,
%$\mathfrak{sp}_{4}$, $\mathfrak{o}_4$ степени многочлена $f$ по
%переменным $z_{-2,-1},z_{-2,1}$ определяются левыми элементами в
%таблице Гельфанда-Цетлина, кодирующих вектор, отвечающий многочелну
%$f$.

Using an explicit description of admissible operators one obtains a proposition

\begin{prop}\label{p4rg}
In the case $g_n$ all polynomial soltions $I_n$  can be written as
$\Omega^{g_{n}}f$, they are encoded by tableaux of type

\begin{align}\begin{split}\label{dcgln2}
&m_{-n,n},\,\,\,\,m_{-n+1,n},\,\,\, m_{-n+2,n},,....\,\,\,\,\,\,\,\,\,m_{-3,n},\,\,\,\,\\
&\,\,\,\,m'_{-n,n},\,\,\, m'_{-n+1,n},,....\,\,\,\,\,\,\,\,\,\,\,\,\, m'_{-3,n},\,\,\,\,  D \\
&\,\,\,\,\,\,\,\,m_{-n,n-1},\,\,\, m_{-n+1,n-1},,....\,\,\,\,\,\,\,\,\,\,\,\,\, m_{-3,n-1},\,\,\,\, \,\,\,\, ,
\end{split}\end{align}
where $D$  is a Gelfand-Tsetlin type  tableaux  for a $g_2$-representation with the highest weight $[m_{-2,n},m_{-1,n}]$. The following inequalities must take place

 \begin{align}\begin{split}\label{ineq}
& m_{-n,n}\geq m'_{-n,n}\geq  m_{-n+1,n}\geq....\geq m_{-3,n}\geq m'_{-3,n}\geq m_{-2,n} \geq...\\
& m'_{-n,n}\geq m_{-n,n-1}\geq  m'_{-n+1,n}\geq .... \geq m'_{-3,n} \geq m_{-3,n-1} \geq m'_{-2,n}\geq ...
% m_{-n,n-2}\geq  m_{-n+1,n-2}\geq ,.... \geq m_{-3,n-2} \geq
% m_{-2,n-2}\geq...
 ,\end{split}
 \end{align}

supplied with the inequalities corresponding to the  $g_2$-tableau $D$.

The polynomial $f$ is defined by $D$,  and $\Omega^{g_{n}}$ is defined by the rest part of the tableau \eqref{dcgln2}.

\end{prop}

\subsection{The Gelfand-Tsetlin type base for $g_{n}$}

 To use this indexation of
 $g_{n-1}$-highest  obtain in Proposition \ref{p4rg} vectors for the construction of the Gelfand-Tsetlin type base we must prove the following Proposition.

\begin{prop}
The lower row in   \eqref{dcgln2} is the $g_{n-1}$-weight of the corresponding
$g_{n-1}$-highest vector.
\end{prop}
\proof The proof is based on the following fact. Let us be given a polynomial from
$\mathfrak{k}[z_{-3,-1},z_{-3,1},...,z_{-n,-1},z_{-n,1}]$,
corresponding to the tableau   \eqref{dcgln2}. The action on it of the operators $E_{i,i}$ $i=-n,...,-3$, in the case $\mathfrak{gl}_{n+1}$,
and of the operators $F_{i,i}$, in the case $g_n$, coincide.

For the operators $E_{i,i}$ this polynomial is an eigenvector with the eigenvalue   $m_{i,n-2}$.  Thus  for the components of the weight with indices $i=-n,...,-3$ the statement is proved. For $i=-2$ the statement follows form the analogous statement for $g_2$.

\endproof

 Applying the standard procedure of construction of the Gelfand-Tsetlin base we obtain the Theorem.

\begin{thm}  In a represention of $g_{n}$ there exist a base called the Gelfand-Tsetlin type base. Its vector are encoded by tableaux \eqref{gcmspo}. The rows $[m]_{k}, [m']_{k},[m]_{k-1}$ of these tableaux are of type \eqref{dcgln2}, where $D$ is a Gelfand-Tsetlin tableau for $g_2$.

The elements of these tableaux satisfy the inequalities \eqref{ineq},
supplied with inequalities corresponding to the  $g_2$-tableau $D$.

\end{thm}

Let us give the formulas for the weight of the vector encoded by a tableau $(m)$. Denote it as
\begin{equation}
\Delta(m)=[\Delta(m)_{-n},...,\Delta(m)_{-1}],
\end{equation}

 \begin{prop}\label{strokan}
\begin{align*}
\Delta(m)_{-k+n-1}=-2\sum_{i}m'_{i,k}+\sum_{i}m_{i,k}+\sum_{i}m_{i,k-1} \text { in the cases } \mathfrak{sp}_{2n},\mathfrak{o}_{2n},\\
\Delta(m)_{-k+n-1}=-2\sum_{i}m'_{i,k}+\sum_{i}m_{i,k}+\sum_{i}m_{i,k-1}-\sigma_{-k}\text { in the cases } \mathfrak{o}_{2n+1}\\
\end{align*}

\end{prop}

The proof can be found in  Appendix in Section
\ref{sootvn}.

\subsection{Reduced matrix elements}

\label{reductgn}

To calculate Wigner coefficients and  matrix elements of generators we need  reduced matrix elements of  the operator  $F_{-1,-2}$,
viewed as $g_{n-1}$-tensor operator that acts between
$g_{n-1}$-representations into which a
$g_n$-representation splits.

Let  $(\bar{m})$, $(m)$be two tableaux for the algebra $g_n$.
In the calculation of the reduced matrix elements we can suggest that these tableaux are maximal with respect to the subalgebra $g_{n-1}$.

\begin{defn}
Denote as $(m)_{red}$ the part of the tableau $(m)$, formed by three upper rows   (that is the tableau \eqref{dcgln2}) from  which the
 $g_2$-tableau $D$ is removed. Thus the three upper rows of  $(m)$ (on which the reduced matrix element depend) can be written as $(m)_{red}D$.
\end{defn}

Note that $(m)_{red}$ is also a part of a
$\mathfrak{gl}_{n+1}$-tableau. Its three upper rows can be written as $(m)_{red}D$, where  $D$ is a $\mathfrak{gl}_3$-tableau. The action of $E_{-1,-2}$ on the tableau
  $(m)_{red}D$  is know (see \cite{zh}). The result is a linear combination of tableaux, the  $i$-th tableau in this combination is obtained by subtracting  $1$ from the $i$-th element of the third row.

Thus one can write

 \begin{equation}\label{expr1} E_{-1,-2}((m)_{red}D)=(E_{-1,-2}(m)_{red})D+(m)_{red}(E_{-1,-2}D).
\end{equation}

To obtain an analogous equality for $F_{-1,-2}$,
let us use the formula for the polynomial
$f$, corresponding to $(m)_{red}D$ that was obtained in \ref{indn}. One has

\begin{equation}\label{expr2} f=\Omega^{g_n}c,
\end{equation}

where $c$is a polynomial, corresponding to $D$.

The action of  $E_{-1,-2}$ on  $f$ is described as follows. We multiply the expression  \eqref{expr2} onto  $E_{-1,-2}$ in the left and then move  $E_{-1,-2}$ to $f_0$. From the commutation relations new summands appear. They correspond to the member
$(E_{-1,-2}(m)_{red})D$ in the expression \eqref{expr1}.
 The  member corresponding to the action of  $E_{-1,-2}$   on  $f_0$, corresponds to the member $(m)_{red}(E_{-1,-2}D)$ in the expression \eqref{expr1}.

% То есть редуцированный матричный элемент зависит только от трёх верхних строк этих таблиц, %обозначим их  $(\bar{m})_{red}$, $(m)_{red}$.

%Пусть $(\bar{m})_{red}$, $(\hat{m})_{red}$ - две таблицы для алгебры $\mathfrak{gl}_{n+1}$,
% совпадающие с  $(\bar{m})_{red}$, $(m)_{red}$. При этом в случае $\mathfrak{o}_{2n+1}$ число $\sigma_{-n+1}$ отброшено.
% Справедлива следующая лемма.

Let us formulate the ruler for calculation of  $<(\bar{m})_{red}|F_{- 1, -2}|(m)_{red}>_{red}$.

\begin{lem}\label{redgn}
\begin{equation} F_{-1,-2}((m)_{red}D)=(E_{-1,-2}(m)_{red})D+(m)_{red}(F_{-1,-2}D).
\end{equation}
\end{lem}

\proof
 The proof is an immediate consequence  of the following fact. The correspondence: $E_{i,-2}\mapsto F_{i,-2}$, $i=-n,...,-1$ is agreed with commutators.

\endproof

 Let us give an explicit formula for the reduced matrix elements.

 Note that when we apply the operator

\begin{equation} 
\Omega^{\mathfrak{gl}_{n+1}}=\prod_{i=-n}^{-3}
\nabla_{-1,i}^{m'_{i,n}-m_{i,n-1}}
  \prod_{i=-n}^{-3}z_{i,1}^{m_{i,n}-m'_{i,n}},
\end{equation} 

to the function, which is identically  equal to one, we get a vector corresponding to the diagram

%\begin{align}\begin{split}
%&m_{-n,n},\,\,\,\,m_{-n+1,n},\,\,\, m_{-n+2,n},,....\,\,\,\,\,\,\,\,\,m_{-3,n},\,\,\,\, 0,\,\,\,\, 0,\,\,\,\, 0\\
%&\,\,\,\,m_{-n,n-1},\,\,\, m_{-n+1,n-1},,....\,\,\,\,\,\,\,\,\,\,\,\,\, m_{-3,n-1},\,\,\,\,0,\, 0,\,\,\,\, \\
%&\,\,\,\,\,\,\,\,m_{-n,n-2},\,\,\,
%m_{-n+1,n-2},,....\,\,\,\,\,\,\,\,\,\,\,\,\, m_{-3,n-2},\,\,\,\,
%\,\,\,\,\,0.\,\,\,\,\,\,\,\,
%\end{split}\end{align}

\begin{align}\begin{split}
&m_{-n,n},\,\,\,\,m_{-n+1,n},\,\,\, m_{-n+2,n},,....\,\,\,\,\,\,\,\,\,m_{-3,n},\,\,\,\, m_{-2,n},\,\,\,\, m_{-1,n},\,\,\,\, m_{1,n}\\
&\,\,\,\,m'_{-n,n},\,\,\, m'_{-n+1,n},,....\,\,\,\,\,\,\,\,\,\,\,\,\, m'_{-3,n},\,\,\,\,m_{-2,n},\, m_{-1,n},\,\,\,\, \\
&\,\,\,\,\,\,\,\,m_{-n,n-2},\,\,\,
m_{-n+1,n-1},,....\,\,\,\,\,\,\,\,\,\,\,\,\, m_{-3,n-1},\,\,\,\,
\,\,\,\,\,m_{-2,n},\,\,\,\,\,\,\,\,
\end{split}\end{align}

where the $\mathfrak{gl}_3$-tableau on the right is maximal.  Denote this $\mathfrak{gl}_3$-tableau as  $D_{max}$.

Define the symbol $<(\bar{m})_{red}  \mid E_{-1,-2} \mid
 (m)_{red}>_{red}$ as follows. If these exists  $i_1=-n,...,-3$,
 such that $\bar{m}_{i_1,n-1}=m_{i_1,n-1}-1$ and all other elements of rows $(\bar{m})_{red}$, $(m)_{red}$ coincide that it equals the  $\mathfrak{gl}_{n+1}$-matrix element $<(\bar{m})_{red}D_{max}\mid E_{-1,-2}\mid (m)_{red}D_{max}>$. Otherwise it is zero.

An explicit formula for the matrix element $<(\bar{m})_{red}D_{max}\mid
E_{-1,-2}\mid (m)_{red}D_{max}>$  is obtained in \cite{zh}. It equals

%\begin{equation}
%\frac{\prod_{i=-n}^{-3} (l_{i,n}-l_{i_1,n-1})}{ \prod_{i=-n,i\neq
%i_1}^{-3} (l_{i,n-1}-l_{i_1,n-1})  },
%\end{equation}

\begin{equation}
\frac{\prod_{i=-n}^{-1} (l'_{i,n}-l_{i_1,n-1})}{ \prod_{i=-n,i\neq
i_1}^{-2} (l_{i,n-1}-l_{i_1,n-1})  },
\end{equation}

where

 \begin{equation}
l'_{i,n}=m'_{i,n}+i,\,\,\,\, l_{i,n-1}=m_{i,n-1}+i.
\end{equation}

The following theorem takes place. It is a direct corollary of Lemma
 \ref{redgn}.

\begin{thm}
\label{t2gn}

\begin{align}\begin{split}&
<(\bar{m})_{red} \bar{D} \mid F_{-1,-2} \mid
(m)_{red}D>_{red}=\delta_{D,\hat{D}}<(\bar{m})_{red}  \mid E_{-1,-2}
\mid (m)_{red}>_{red}+\\&+\delta_{(\bar{m})_{red},(m)_{red}}<\bar{D}
\mid F_{-1,-2} \mid D>_{red}.
\end{split}
\end{align}

It is suggested that
$\bar{m}_{i_1,n-1}=m_{i_1,n-1}-1$ and $\bar{m}_{j,n-1}=m_{j,n-1}$ for
$j\neq i$.

\end{thm}

\subsection{Wigner coefficients}
\label{coew}

Let us obtain formulas for the Wigner coefficients for $g_{n-1}$.
 Following \cite{1963}, let us first obtain formulas for the coefficients

\begin{align}
\label{vig} <(\bar{m}) \begin{pmatrix} j \\ [1\dot{0}]_{n-1} \\ -2
\end{pmatrix} (m)>.
\end{align}
% где $i<-1$.

We use the following fact.
 Let us be given a diagram $(m)$ that define a vector in s  $g_{n-1}$-representation
  with the highest weight $[m_{-n,n-1},...,m_{-2,n-1}]$.

A polynomial on the group $Z_{n-1}$ corresponds to this vector. Consider it as a polynomial on a bigger group  $Z_n$.

\begin{prop}\label{predl}
The vector that corresponds to this polynomial belong to a  $g_{n}$-representation with the highest weight $[m_{-n,n-1},...,m_{-2,n-1},0]$. The corresponding tableau is of type $\begin{pmatrix}max\\ m\end{pmatrix}$, this is a  $g_{n}$-tableau that is obtained form
    $(m)$ by adding two maximum row.

\end{prop}

The proof can be found in Appendix in Section
\ref{docpredl}.

Now return to the calculation of the Wigner coefficient
 \eqref{vig}. Let us be given two $g_{n-1}$-tableaux $(\bar{m})$ and  $(m)$. We can suggest that they are $g_{n-2}$-maximal. The three upper rows of these tableau are of type $(\bar{m})_{red}\bar{D}$ and $(m)_{red}D$.  As in Proposition \ref{predl}
   add to them two maximal
  $g_n$-rows, denote the five rows that we obtain as $\widetilde{(\bar{m})_{red}\bar{D}}$ and $\widetilde{(m)_{red}D}$.

  Let us prove the equality.

\begin{prop}\label{pr12}
\begin{equation}
 <(\bar{m}) \begin{pmatrix} j \\ [1\dot{0}]_{n-1} \\ -2
\end{pmatrix} (m)>=<\widetilde{((\bar{m})_{red}\bar{D})} \mid F_{-1,-2} \mid  \widetilde{((m)_{red}D)}>.
\end{equation}
\end{prop}

\proof
Let us apply the Wigner-Eckart theorem to the matrix element on the right. It equals to the product of a reduced matrix element
 and a Wigner coefficient that occurs on the left side of the equality. Since the upper two rows of tableaux $(\bar{m})_{red}\bar{D})$, $\widetilde{((m)_{red}D)}$ are maximal the reduced matrix element equals to $1$. This proves the equality.

\endproof

  Let us calculate the matrix element that occurs on the right in the equality \ref{pr12}.

Let us be given a tableau that defines a
$g_{n-2}$-highest vector in a $g_n$-representation. This is a tableau of type

\begin{align}
\begin{split}\label{diagc}
& m_{-n,n},\,\,\,....\,\,\,m_{-4,n}\,\,\,  \\
&\,\,\,m'_{-n,n},\,\,\, ....\,\,\,m'_{-4,n}\,\,\,\,\\
&\,\,\,\,\,\,m_{-n,n-1},\,\,\, ....\,\,\,m_{-4,n-1}\,\,\,\, C\\
&\,\,\,\,\,\,\,\,\,m'_{-n,n-1},\,\,\, ....\,\,\,m'_{-4,n-1}\,\,\,\,\\
&\,\,\,\,\,\,\,\,\,\,\,\,m_{-n,n-2},\,\,\, ....\,\,\,m_{-4,n-2}\,\,\,\,,\\
\end{split}
\end{align}

 where $C$ is a $g_3$-tableau. Denote the tableau \eqref{diagc} shortly as
$(k)C$.
 Using the technique of raising operators $\nabla_{ij}$ (see \cite{zh}),  and applying the  arguments that were used in the poorf of the formula \eqref{f52}, one obtains that  the polynomial that corresponds to this tableau
  can be written as
$\Omega f,$ where $f$ is a polynomial that corresponds to a
$g_3$-tableau $C$ and the operator $\Omega $  is defined as follows

\begin{equation}
\Omega =\prod_{i=-n}^{-3} \nabla_{-3,i}^{m'_{i,n-1}-m_{i,n-2}}
\prod_{i=-n}^{-3}
\nabla_{-2,i}^{m_{i,n-1}-m'_{i,n-1}}\prod_{i=-n}^{-3}
\nabla_{-1,i}^{m'_{i,n}-m_{i,n-1}}
  \prod_{i=-n}^{-3}z_{i,1}^{m_{i,n}-m'_{i,n}}
\end{equation}

The considered tableau $\widetilde{((m)_{red}D)}$ is of type
\eqref{diagc}. One has $(n)=\widetilde{(m)_{red}}$ and
$C=\widetilde{D}$. By analogy with the proof of Lemma \ref{redgn} one concludes that

\begin{equation}
 F_{-1,-2}\widetilde{((m)_{red}D)}=\widetilde{(E_{-1,-2}(m)_{red})}\widetilde{D}+\widetilde{(m)_{red}}F_{-1,-2}\widetilde{D}.
\end{equation}

 When one passes to matrix elements, one gets

 \begin{align}\begin{split}&
 <\widetilde{((\bar{m})_{red}\bar{D})} \mid F_{-1,-2} \mid  \widetilde{((m)_{red}D)}>=
 \delta_{\bar{D},D}<\widetilde{((\bar{m})_{red})} \mid E_{-1,-2} \mid  \widetilde{((m)_{red}}>+\\
 &+\delta_{(\bar{m})_{red},(m)_{red}}< \widetilde{\bar{D}} \mid F_{-1,-2} \mid \widetilde{D}>.
 \end{split}
 \end{align}

 Thus we have proved a Theorem

\begin{thm}
\label{t3gn}
\begin{align}
\begin{split}
\label{vign}& <(\bar{m}) \begin{pmatrix} j \\ [1\dot{0}]_{n-1} \\ -2
\end{pmatrix} (m)>=\\&=
\delta_{\bar{D},D}<(\bar{m})_{red} \mid E_{-1,-2} \mid  (m)_{red}>+\delta_{(\bar{m})_{red},(m)_{red}}<\bar{D} \mid F_{-1,-2} \mid D>.
\end{split}
\end{align}

\end{thm}

 Let us give rulers for calculation of summands that occur on the right hand side in Theorem  \ref{t3gn}.

\subsubsection{The matrix element $<\widetilde{((\bar{m})_{red})} \mid E_{-1,-2} \mid  \widetilde{((m)_{red}}>$}

As in previous section this matrix element can expressed through a matrix element of the algebra $\mathfrak{gl}_{n+1}$.

\begin{thm} If there exists
$i_1=-n,...,-3$,
such that $\bar{m}_{i_1,n-1}=m_{i_1,n-1}-1$, and all other elements of rows $(\bar{m})_{red}$, $(m)_{red}$ coincide than the considered matrix
elements equals

\begin{equation}
\frac{\prod_{i=-n}^{-1} (l_{i,n-1}-l_{i_1,n-2})}{ \prod_{i=-n,i\neq
i_1}^{-2} (l_{i,n-2}-l_{i_1,n-2})  },
\end{equation}

where

 \begin{equation}
l_{i,n-1}=m_{i,n-1}+i,\,\,\,l_{i,n-2}=m'_{i,n-1}+i.
\end{equation}

If thexe exist no such index than the considered matrix element equals zero

\end{thm}

  \subsubsection{The matrix element $< \widetilde{\bar{D}} \mid F_{-1,-2} \mid \widetilde{D}>$. The case $\mathfrak{o}_5$}

This matrix element equals to a $\mathfrak{o}_5$-Wigner coefficiens

\begin{equation}
< \widetilde{\bar{D}} \mid F_{-1,-2} \mid \widetilde{D}>=<\bar{D} \begin{pmatrix} j \\ [10]_{\mathfrak{o}_5} \\ -2
\end{pmatrix}  D>.
\end{equation}

Let us calculate it directly.  One can suggest that
$\mathfrak{o}_5$-tableaux $\bar{D}$ and $D$ are
$\mathfrak{o}_3$-maximal.

To the tableau $D$ there corresponds a polynomial  $f$ on $Z_{\mathfrak{o}_5}$.
ТSince $D$ is $\mathfrak{o}_3$-maximal, then
$f=f(z_{-3,-2},z_{-3,2})$. By Proposition \ref{predl},
to the tableau $\widetilde{D}$ there corresponds the same polynomial, but considered as a polynomial on  $Z_{\mathfrak{o}_5}$.

Let find the action of the operator $e^{tF_{-1,-2}}$ on the polynomial $f$.
The explicit calculation gives that

\begin{equation} \label{o5wig}
f(z_{-3,-2},z_{-3,2})\mapsto
(1+tz_{-2,-1})^{m_{-2,2}-m_{-1,2}}(z_{-3,-2}+tz_{-3,-1},z_{-3,2})
\end{equation}

When we were defining the Gelfand-Tsetlin type base for
$\mathfrak{o}_5$ in Section \ref{indb} for the polynomial $f$ we have constructed a polynomial  $f^*$ on the group $Z_{\mathfrak{gl}_4}$.

 By explicit calculations it can be shown that on the polynomial   $f^*$ the operator   $e^{tE_{-1,-2}}$ acts by the same formula \eqref{o5wig}.  Thus the correspondence conjugates  the actions of  $e^{tF_{-1,-2}}$ and $e^{tE_{-1,-2}}$.

In Section \ref{gcbaseo5} 
using this correspondence the Gelfand-Tsetlin type base for $\mathfrak{o}_5$ was constructed. To a  $\mathfrak{o}_5$-tableau $D$, that defines a $\mathfrak{o}_3$-highest vector (and a polynomial $f$), there corresponds a  $\mathfrak{gl}_3$-tabelau $D^*$,  which is obtained from $D$ by removing the zero from the upper row and $\sigma_{-2}$ (to the tableau $D^*$ there corresponds the polynomial $f^*$).

Thus we have

\begin{equation}\label{o5gl3wig}
< \widetilde{\bar{D}} \mid F_{-1,-2} \mid \widetilde{D}>=< \widetilde{\bar{D^*}} \mid E_{-1,-2} \mid \widetilde{D^*}>
\end{equation}

The matrix element on the right in \eqref{o5gl3wig} is a  $\mathfrak{gl}_3$-Wigner coefficient.

Thus we have proved the theorem
 \begin{thm}

The Wigner coefficient for $\mathfrak{o}_5$ and $\mathfrak{gl}_3$ are equal

\begin{equation}
<\bar{D} \begin{pmatrix} j \\ [10]_{\mathfrak{o}_5} \\ -2
\end{pmatrix}  D>=<\bar{D^*} \begin{pmatrix} j \\ [100]_{\mathfrak{gl}_3} \\ -2
\end{pmatrix} D^*>=\begin{vmatrix} j: 3 \\  2\end{vmatrix}^{[m]_2,[m']_{2}},
\end{equation}
where  $[m]_2=[m_{-2,2},m_{-1,2},0]$  и $[m']_2=[m'_{-2,2},m'_{-1,2}]$.

\end{thm}

\begin{cor}

\begin{equation}
< \widetilde{\bar{D}} \mid F_{-1,-2} \mid \widetilde{D}>=\begin{vmatrix} j: 3 \\  2\end{vmatrix}^{[m]_2,[m']_{2}},
\end{equation}
where  $[m]_2=[m_{-2,2},m_{-1,2},0]$  и $[m']_2=[m'_{-2,2},m'_{-1,2}]$.
\end{cor}

  \subsubsection{The matrix element $< \widetilde{\bar{D}} \mid F_{-1,-2} \mid \widetilde{D}>$. The case $\mathfrak{sp}_4$}

In this case the matrix element equals to a $\mathfrak{sp}_4$-Wigner coefficient

\begin{equation}
< \widetilde{\bar{D}} \mid F_{-1,-2} \mid \widetilde{D}>=<\bar{D} \begin{pmatrix} j \\ [10]_{\mathfrak{sp}_4} \\ -2
\end{pmatrix}  D>.
\end{equation}

In section \ref{gcbasesp4} the Gelfand-Tsetlin type base for  $\mathfrak{sp}_4$ was defined. To a  $\mathfrak{sp}_4$-tableau $D$
 there corresponds a  $\mathfrak{gl}_3$-tableau $D^*$ from which the zero in the upper row is removed.

Analogously to the case $\mathfrak{o}_5$ one can prove

 \begin{thm}

The Wigner coefficient for $\mathfrak{sp}_4$ and $\mathfrak{gl}_3$ are equal

\begin{equation}
<\bar{D} \begin{pmatrix} j \\ [10]_{\mathfrak{sp}_4} \\ -2
\end{pmatrix}  D>=<\bar{D^*} \begin{pmatrix} j \\ [100]_{\mathfrak{gl}_3} \\ -2
\end{pmatrix} D^*>=\begin{vmatrix} j: 3 \\  2\end{vmatrix}^{[m]_2,[m']_{2}},
\end{equation}

where $[m]_2=[m_{-2,2},m_{-1,2},0]$  и $[m']_2=[m'_{-2,2},m'_{-1,2}]$.

\end{thm}

\begin{cor}

\begin{equation}
< \widetilde{\bar{D}} \mid F_{-1,-2} \mid \widetilde{D}>=\begin{vmatrix} j: 3 \\  2\end{vmatrix}^{[m]_2,[m']_{2}},
\end{equation}
where $[m]_2=[m_{-2,2},m_{-1,2},0]$  и $[m']_2=[m'_{-2,2},m'_{-1,2}]$.
\end{cor}

\subsubsection{The matrix element $< \widetilde{\bar{D}} \mid F_{-1,-2} \mid \widetilde{D}>$. The case $\mathfrak{o}_4$}

In this case the matrix element equals to a  $\mathfrak{o}_4$-Wigner coefficient.

\begin{equation}
< \widetilde{\bar{D}} \mid F_{-1,-2} \mid \widetilde{D}>=<\bar{D} \begin{pmatrix} j \\ [10]_{\mathfrak{o}_4} \\ -2
\end{pmatrix}  D>.
\end{equation}

Let us calculate directly the $\mathfrak{o}_4$-Wigner coefficient on the right.

The index $j$ can take values $-2$  and $-1$. If $j=-2$ then

\begin{equation}
\bar{m}_{-2}=m_{-2}+1,\,\,\,\bar{m}_{-1}=m_{-1},\\
\end{equation}

  and if $j=-1$ then

\begin{equation}
\bar{m}_{-2}=m_{-2},\,\,\,\bar{m}_{-1}=m_{-1}+1.\\
\end{equation}

For $r_{-2}=m_{-2}-m_{-1}$ and $r_{-1}=m_{-2}+m_{-1}$ (these  are
highest weights for two  $\mathfrak{sl}_2$ copies), one has in the
case $j=-2$

\begin{equation}
\bar{r}_{-2}=r_{-2}+1,\,\,\,\bar{r}_{-1}=r_{-1}+1,\\
\end{equation}

and in the case $j=-1$

\begin{equation}
\bar{r}_{-2}=r_{-2}-1,\,\,\,\bar{r}_{-1}=r_{-1}-1.\\
\end{equation}

The weights $p$  and $q$ for two copies of $\mathfrak{sl}_2$ are expressed through the elements of a  $\mathfrak{o}_4$-tableau using the Proposition \ref{pq}.

Thus one gets the Theorem

\begin{thm}

Put $[m^1]_2=[\frac{m_{-2,2}-m_{-1,2}}{2},0]$,     $[m^2]_2=[\frac{m_{-2,2}+m_{-1,2}}{2},0]$.

For $j=-2$ one has

\begin{equation}
<\bar{D} \begin{pmatrix} j \\ [10]_{\mathfrak{o}_4} \\ -2
\end{pmatrix}  D>=\begin{vmatrix} \frac{1}{2}: 2 \\  1\end{vmatrix}^{[m^1]_2, p}\begin{vmatrix} \frac{1}{2}: 2 \\  1 \end{vmatrix}^{[m^2]_2, q}.
\end{equation}

For $j=-1$ one has

\begin{equation}
<\bar{D} \begin{pmatrix} j \\ [10]_{\mathfrak{o}_4} \\ -2
\end{pmatrix}  D>=\begin{vmatrix} \frac{1}{2}: 2 \\  1 \end{vmatrix}^{[m^1]_2, p}\begin{vmatrix} -\frac{1}{2}: 2 \\  1\end{vmatrix}^{[m^2]_2, q}.
\end{equation}

\end{thm}

\subsection{Reduced Wigner coefficients}

For the Wigner coefficients \begin{align} \label{vigui}
<(\bar{m})
\begin{pmatrix} j \\ [1\dot{0}]_{n-1} \\ i
\end{pmatrix} (m)>
\end{align}

the following formula takes place

\begin{align} \label{viguii} <(\bar{m})
\begin{pmatrix} j_1 \\ [1\dot{0}]_{n-1} \\ i
\end{pmatrix} (m)>=\begin{pmatrix} j_{-i-1}  \\ [1\dot{0}]_{n-i} \\  i \end{pmatrix}\prod_{l=1}^{-i-2}\begin{vmatrix} j_l  \\ [1\dot{0}]_{n-l}  \\ j_{l+1}
\end{vmatrix}
\end{align}

To calculate all Wigner coefficients we must obtain a formula for the reduced Wigner coefficients $\begin{vmatrix} j_l  \\ [1\dot{0}]_{n-l}  \\ j_{l+1}
\end{vmatrix} $.
%участвующих в формуле \eqref{viguii}.

 Let s obtain the formula for the reduced Wigner coefficients using the previous calculations. The following theorem takes place.

\begin{thm}
\label{t4}
\begin{align}
\begin{vmatrix}j_1 \\ [1\dot{0}]_{n-1}  \\ j_2 \end{vmatrix}=\begin{pmatrix}j_1 \\ [1\dot{0}]_{n-1}  \\ -2 \end{pmatrix} <(\bar{m'})_{red} \bar{D}' \mid F_{-2-3} \mid
(m')_{red}D'>_{red},
\end{align}
where $'$ means that we take only the part of the tableaux that correspond to $g_{n-1}$.
\end{thm}
\proof

To prove the theorem let us calculate the matrix element $<(\bar{m})
\mid F_{-1-3} \mid (m)>$, $i<-2$ in two ways.

Firstly apply to the matrix element the Wigner-Eckart theorem and decompose the Wigner coeficient into the product of a reduced Wigner coefficient and a Wigner coefficient.
One has

\begin{align}
<(\bar{m}) \mid F_{-1-3} \mid (m)>=<(\bar{m})_{red}\bar{D} \mid F_{-1-3} \mid (m)_{red}D>_{red}
\begin{vmatrix}j_1 \\ [1\dot{0}]_{n-1}  \\ j_2 \end{vmatrix}
 \begin{pmatrix}j_2 \\ [1\dot{0}]_{n-2}  \\ -3 \end{pmatrix}.
\end{align}

Secondly using the commutation relation
$F_{-1,-3}=[F_{-1,-2},F_{-2,-3}]$  we obtain
\begin{align}\begin{split}&
<(m') \mid F_{-1,-3} \mid (m)>=<(\bar{m})_{red} \bar{D}\mid F_{-1-2} \mid
(m)_{red}D>_{red}
\begin{pmatrix}j_1 \\ [1\dot{0}]_{n-1}  \\ -2 \end{pmatrix} \cdot \\& \cdot<(\bar{m'})_{red}\bar{D}' \mid F_{-2-3} \mid
(m')_{red}D'>_{red}
\begin{pmatrix}j_2 \\ [1\dot{0}]_{n-2}  \\ -3
\end{pmatrix},\end{split}
\end{align}

where $'$ means that we take only the part of the tableaux that correspond to $g_{n-1}$.

Compare two expressions, one obtains

\begin{align}
\begin{vmatrix}j_1 \\ [1\dot{0}]_{n-1}  \\ j_2 \end{vmatrix}=\begin{pmatrix}j_1 \\ [1\dot{0}]_{n-1}  \\ -2 \end{pmatrix} <(\bar{m'})_{red} \bar{D}' \mid F_{-2-3} \mid
(m')_{red}D'>_{red}
\end{align}

The theorem is proved

 \endproof

The expressions in Theorem \ref{t4} are obtained in previous Sections.

\subsection{Matrix elements}

%В разделе \ref{singsp2n} было дано построение базиса типа
%Гельфанда-Цетлина.

Using the previous results let us write the formulas for the action of generators of  $g_n$ in the base that we have constructed. It is enough to give a formula for the action of  $F_{-1,-2}$.
%Они легко получаются с помощью формул
%Вигнера-Эккарта.
The following theorem takes place.

\begin{thm}

\label{matgn}

\begin{align}
\begin{split} &
<(\bar{m})| F_{-1,-2} |(m)>=<(\bar{m})_{red} \bar{D} \mid F_{-1,-2} \mid
(m)_{red}D>_{red}<(\bar{m}) \begin{pmatrix} j \\ [1\dot{0}]_{n-1} \\ -2
\end{pmatrix} (m)>,
\end{split}
\end{align}

where the expression for the factors are given in theorems \ref{t2gn} and
\ref{t3gn}.

\end{thm}

\proof The theorem is proved by application of the Wigner-Eckart theorem

\endproof

\section{Appendix}
\label{appendi}

\subsection{Lie algebras.}

\label{alglie}

%Напомним определения используемых алгебр ли $\mathfrak{sp}_{2n}$,
%$\mathfrak{o}_N$, см также \cite{jb}.

\subsubsection{The symplectic algebra $\mathfrak{sp}_{2n}$}

 Take the space $\mathbb{C}^{2n}$. Choose a base and let us index its elements by numbers
 $-n,...,-1,1,...,n$.  Fix a skew-symmetric form $$\omega=\sum_{i=1}^n x_i\wedge x_{-i}.$$

 The group $Sp_{2n}$ consists of isomorphisms of   $\mathbb{C}^N$, that preserve this skew-symmetric form.  It's Lie algebra is denoted as
$\mathfrak{sp}_{2n}$.

Define the generators of the Lie algebra
$$F_{ij}=E_{ij}-sign (i) sign (j) E_{-j-i}.$$

The only relations between them are
$$F_{ij}=-sign(i)sign(j)F_{-j-i}.$$

The generators $F_{ij} $,
in the case $i>j$,
correspond to negative roots. The generators $F_{ij} $,
in the case $i<j$,
correspond to positive roots. In the case $i=j$ the generator belongs to the Cartan subalgebra.

\subsubsection{The orthogonal algebra $\mathfrak{o}_{N}$}

Take the space $\mathbb{C}^N$.  Let $n$ be such that $N=2n$ in the case of even  $N$, and $N=2n+1$ in the case of odd  $N.$

Choose a base in $\mathbb{C}^N$ and index its elements by  $-n,...,-1,1,...,n$ in the case of even  $N$ and by numbers
$-n,...,-1,0,1,...,n$    in the case of odd $N$.

 In the case of even $N$ take a quadratic form  $$x_{-n}x_{n}+...+x_{-1}x_1,$$ в
  and in the case of odd  $N$ take a quadratic form $$x_{-n}x_n+...+x_{-1}x_1+x_0^2.$$

The group $O_N$ consists of isomorphisms of  $\mathbb{C}^N$,
that preserve this quadratic form. It's Lie algebra is denoted as
$\mathfrak{o}_{N}$.

Introduce  generators
$$F_{ij}=E_{ij}- E_{-j-i}.$$

The only relations between them are
$$F_{ij}=-F_{-j-i}.$$

The generators $F_{ij} $,
in the case $i>j$,
correspond to negative roots. The generators $F_{ij} $,
in the case $i<j$,
correspond to positive roots. In the case $i=j$ the generator belongs to the Cartan subalgebra.

%Иными словами, если $G$  есть матрица квадратичной формы, то матрица
%$M\in O_N$ тогда и только тогда, когда выполняется равенство

%$$M^tGM=G.$$

%При такой реализации ортогональной группы, нижнетреугольные матрицы
%$M\in O_N$ являются повышающими операторами, верхнетреугольные - понижающими. Диагональные матрицы  принадлежам тору.

\subsection{Tensor operators and Wigner coefficients }

In this section the Wigner coefficients are defined, the solution of the multiplicity problem is given, the reduced Wigner coefficients are introduced.  We follow the analogous discussion for the case of the algebra
 $\mathfrak{gl}_n$ in \cite{1968}.

The Wigner coefficients are closely related with irreducible tensor operators.
\begin{defn}
An irreducible tensor operator of type $[M]_n$,where $[M]_n$ is a  dominant $g_n$-weight, is an indexed by vectors $(M)\in
V^{[M]_n}$the set of linear mappings

$$f_{(M)}: V^{[\bar{m}]_{n}}  \rightarrow V^{[m]_{n}},$$

which has the following property. For $g\in g_n$ one has
 $$[g,f_{(M)}]=f_{g(M)}.$$

\end{defn}

For given  $[\bar{m}]_n$, $[m]_n$ and $[M]_n$ the tensor operator $ V^{[\bar{m}]_{n}}
\rightarrow V^{[m]_{n}} $ of type $[M]_n$ is not unique. If it is not unique one says that for this tensor operator the multiplicity problem takes place.

By Wigner-Eckart theorem the matrix elements of a tensor operator decompose into a product of a factor that depents only on highest weights
весов $[m]_{n},[\bar{m}]_{n},[M]_{n}$ (this factor is called the reduced matrix element) and a Wigner coefficient that defines an interwinnig operator $\Phi:
V^{[\bar{m}]_{n}}\rightarrow V^{[M]_{n}}\otimes V^{[m]_{n}}$.

\subsection{ The solution of the multiplicity problem}
As it is know the interwinnig operator  $\Phi: V^{[\bar{m}]_{n}}
\rightarrow V^{[M]_{n}}\otimes V^{[m]_{n}}$
 is in general not unique. On the other language this means that one irreducible representation $V^{[\bar{m}]_{n}} $ can occur in splitting of tensor product not once. But the parametrization of all interwinnig operators is well-known.

All such operators are indexed by tableaux $(\Gamma)\in V^{[M]_{n}}$,
such that $$[\bar{m}]_n=\Delta (\Gamma)+[m]_n,$$ where $\Delta (\Gamma)$ is the weight of  $v$. This vector   $\Phi$ is defined by the operator as follows

$$\Phi((m)_{max})=\Delta (\Gamma)\otimes (m)_{max}+l.o.t.,$$  where $l.o.t.$ (lower order terms)denotes a sum of tensor products of weight vectors where the second vector has a weight lower than $[m]_n$.
 For the group $U(n)$ this was proved by Biedenharn and Baird in \cite{19633}.

The corresponding Wigner coefficient is denoted as

\begin{equation}\label{viggen}
<\begin{pmatrix}  [\bar{m}]_n \\  (m')_{n-1} \end{pmatrix}
\begin{pmatrix}  (\Gamma)_{n-1} \\ [ M]_n\\ (M)_{n-1}\end{pmatrix}
\begin{pmatrix}[m]_n\\ (m)_{n-1}\end{pmatrix}>
\end{equation}

This coefficient can be nonzero only if the following equality holds $[\bar{m}]_n=\Delta(\Gamma)+[m]_n$.

%Коэффициент Вигнера, как известно, задает неприводимый тензорный оператор

%$$ V^{[m']_{n}} \rightarrow V^{[m]_{n}}.$$

\subsection{Reduced Wigner coefficients}

Take a Wigner coefficient\eqref{viggen} for the algebra  $g_n$.
It also defines a tensor operator for the algebra $g_{n-1}$.
Decompose it into a sum of irreducible tensor operators and apply the Wigner-Eckart theorem. One gets

 \begin{align*}
&<\begin{pmatrix}  [\bar{m}]_n \\  [\bar{m}']{n} \\ [\bar{m}]_{n}\\
(\bar{m})_{n-2}
\end{pmatrix} \begin{pmatrix} ( \Gamma )_{n-1}\\ [ M]_n\\ (M)_{n-1}
\end{pmatrix} \begin{pmatrix}[m]_n\\ [m']_{n}  \\  [m]_n \\ (m)_{n-1}
\end{pmatrix}>=\sum_{ (\gamma)_{n-2}}  <\begin{pmatrix}  [\bar{m}]_n \\  [\bar{m}']_{n} \\ [\bar{m}]_{n}\\ (\bar{m})_{n-2} \end{pmatrix} \begin{vmatrix}  (\Gamma)_{n-1} \\ [ M]_n\\ (\gamma)_{n-1}
\end{vmatrix} \begin{pmatrix}[m]_n\\ [m']_{n}  \\  [m]_{n-1} \\ (m)_{n-2}
\end{pmatrix}> \\&
  <\begin{pmatrix}  [\bar{m}]_{n-1} \\  (\bar{m})_{n-2} \end{pmatrix} \begin{pmatrix}  (\gamma)_{n-2} \\ [ M]_{n-1}\\ (M)_{n-2}
\end{pmatrix} \begin{pmatrix}[ [m]_{n-1} \\ (m)_{n-1}
\end{pmatrix}>
  \end{align*}

In this formula the first factor on the right is a notation for the reduced  matrix element of an irreducible $g_{n-1}$-tensor operator.

This factor is called the reduced Wigner coefficient.

 The tableau $(\gamma)_{n-1}$ is obtained by adding to the tableau $(\gamma)_{n-2}$ the rows  $[M']_{n}$ and $[M]_{n-1}$.
Note that the reduced matrix element does not depend on the rows $[\bar{m}]_{n-2}$, $[m]_{n-2}$ and below.
% Действительно, этот
%коэффициент по очевидным соображениям инвариантен по отношению к
%действию $g_{n-2}$. Кроме того, он инвариантен и по отношению к
%действию алгебры Микельсона-Желобенка $Z(g_{n-1},g_{n-2})$ \cite{M},
%действующей на пространстве старших $g_{n-2}$ векторов.
%Следовательно, коэффициент не зависит также и от
%$[\bar{m}]_{n-2}^{n-1}$ и $[m]_{n-2}^{n-1}$.
 Thus the reduced Wigner coefficient can be denoted as

 \begin{align}\label{redviggen}
<\begin{pmatrix}  [\bar{m}]_n \\  [\bar{m}']{n} \\
[\bar{m}]_{n-1}
\end{pmatrix} \begin{vmatrix}  (\Gamma)_{n-1} \\ [ M]_n\\ (\gamma)_{n-1}
\end{vmatrix} \begin{pmatrix}[m]_n\\ [m']_{n} \\  [m]_{n-1}
\end{pmatrix}>
  \end{align}

This coefficient is nonzero only if the following holds: $[\bar{m}]_n=\Delta(\Gamma)+[m]_n$ и
$[\bar{m}]_{n-1}=\Delta((\gamma)_{n-1})+[m]_n$.

  \subsection{ Fundamental operators}

% В данном
%разделе мы введён некоторый класс тензорных операторов, или
%коэффициентов Вигнера. Для этого класса отсутствует проблема
%кратности.

In the present paper only the Wigner coefficients are considered for which  $[M]_n=[1,0,...,0]=[1 \dot{0}]_{n}$,  that is when   the tensor factor
$V^{[M]_n}$ is a standard representation. Such Wigner coefficients are called fundamental.

Note that weight vectors $(m)$  of the standard  representations are completely defined by their weights $\Delta(m)=[0,...,\pm 1,...,0]$, where $\pm 1$ occurs at the place $i$. If $i=0$ then only $1$ is allowed. Thus the Wigner coefficient is of type can be denoted as

\begin{equation}
\begin{pmatrix} i \\ [1 \dot{0}]_{n}\\j  \end{pmatrix}
\end{equation}

\subsection{The proof of Propositions \ref{stroka05} and \ref{weightofo5}} \label{sootv}

The Proposition \ref{stroka05} is the following statement.

\begin{prop*}
The number $m_{-2,1}$ in the tableau \eqref{gco5} is the
$\mathfrak{o}_3$-weight of the corresponding vector.
\end{prop*}

\proof Take a polynomial
\begin{equation} \label{monomo5}
 \sum_{k,l,r,s}с_{ k,l,r,s } z_{-2,-1}^{k}z_{-2,1}^{l}z_{-1,1}^r z_{-1,0}^s.\end{equation}

Suggest it defines a $\mathfrak{o}_{3}$-highest vector in a
$\mathfrak{o}_{5}$-representaion with the highest  weight $[m_{-2},m_{-1}]$, then to each it's monomial there corresponds a function

\begin{align}\begin{split}
\label{al2}
\delta_{-2}^{m_{-2,2}-k-l}\delta_{-1}^{m_{-1,2}+k-2r-s-l}.\end{split}\end{align}

Then the $\mathfrak{o}_{3}$-weight of the corresponding vector equals
$m_{-2}-k-l.$

Also to \eqref{monomo5} there corresponds a polynomial

\begin{equation} \label{monomgl3}
 \sum_{k,l,r,s}с_{ k,l,r,s }z_{-2,-1}^{k}z_{-2,1}^{l}z_{-1,1}^r
z_{-1,1}^{[\frac{s}{2}]}.\end{equation}

that defines a  $\mathfrak{gl}_{1}$-highest vector in s
$\mathfrak{gl}_{3}$-representaion with the highest weight
$[m_{-2,2},m_{-1,2},0]$, to each it's monomial there corresponds a function

\begin{align}\begin{split}\label{al11} \delta_{-2}^{m_{-2,2}-k-l}\delta_{-1}^{m_{-1,2}+k-r-s}\delta_{1}^{r+l}.\end{split}\end{align}

The weight of the vector that corresponds to this polynomial equals
$m_{-2,2}-k-l.$

%C [m_{-n}-k_{-n}-l_{-n},...,m_{-2}-k_{-2}-l_{-2}]. \end{equation}

Thus the weight of the $\mathfrak{o}_{3}$-highest vector that corresponds to
\eqref{monomo5} equals the weight of the
$\mathfrak{gl}_{1}$-highest vector that corresponds to \eqref{monomgl3}. But the weight of the last vector equals $m_{-2,1}$. The Proposition is proved.

\endproof

Let us prove the Proposition \ref{weightofo5}.

 \begin{prop*}

\label{weighto5}

$\Delta(m)_{-1}=-2\sum_i m'_{i,2}+\sum_i
m_{i,2}+m_{-2,1}-\sigma_{-2}$

\end{prop*}
\proof

Take a polynomial \eqref{monomo5} that corresponds to a
Gelfand-Tsetlin tableaux. Consider the expression \eqref{al2}. One
obtains

\begin{equation}\label{vesochik}
\Delta(m)_{-1}=m_{-1,2}+k-2r-s-l.
\end{equation}

Consider the polynomial the corresponding \eqref{monomgl3} that defines a
$\mathfrak{gl}_1$-highest vector in a
$\mathfrak{gl}_{3}$-representation, its components with indices
$-1$ and $1$  are equal to $m_{-1,2}+k-r-s$  and $r+l$.
Their difference equals $$m_{-1,2}+k-2[\frac{s}{2}]-2-l.$$ This equals
\eqref{vesochik}, if $s$ is even and differs by one from
\eqref{vesochik}, if $s$ is odd. Thus the expression

\begin{equation}
m_{-1,2}+k-2[\frac{s}{2}]-2-l-\sigma_{-1}
\end{equation}
 equals  \eqref{vesochik}.

 Thus $\Delta(m)_{-1}$ is a difference of $\mathfrak{gl}_3$-weights with numbers $-1$ and $1$ minus $\sigma_{-1}$.

%  Чтобы доказать предложение заметить, что  одночлен, соответствующий таблице типа %Гельфанда-Цетлина есть линейная комбинация
 % одночленов
 %  \eqref{monomo5} с одним и тем же весом.

  Note that in the case $\mathfrak{gl}_{3}$ the component of the weight with the index    $-1$  equals $\sum_{i} m'_{i,2}-\sum_i m_{i,2}$, and the component of the weight with the index     $1$  equals
    $m_{-2,1}-\sum_{i} m'_{i,2}$.

\endproof

\subsection{The proof of Proposition \ref{strokan}}

\label{sootvn}

Let us prove the following statement

 \begin{prop*}
\begin{align*}
\Delta(m)_{-k+n-1}=-2\sum_{i}m'_{i,k}+\sum_{i}m_{i,k}+\sum_{i}m_{i,k-1} \text { in the cases } \mathfrak{sp}_{2n},\mathfrak{o}_{2n},\\
\Delta(m)_{-k+n-1}=-2\sum_{i}m'_{i,k}+\sum_{i}m_{i,k}+\sum_{i}m_{i,k-1}-\sigma_{-k}\text
{ in the cases } \mathfrak{o}_{2n+1}
\end{align*}

\end{prop*}
\proof

It is enough to prove the formula for  $\Delta(m)_{-1}$.  One can suggest that the tableau $(\Gamma)$ is maximal with respect to $g_{n-1}$. To this tableau there corresponds a polynomial $f$, the procedure of it's construction is described in Section \ref{indn}.
 The polynomial $f$ is of type $$f=cf_0,$$ where $c\in \mathfrak{k}$ (see Definition of the field $\mathfrak{k}$
 in Section \ref{indn}), and  $f_0\in\mathbb{C}(z_{-3,-1},z_{-3,1},...,z_{-n,-1},z_{-n,1})$.

 The weight of the vector corresponding to   $f$ is calculated as follows (see \cite{zh}).
  To each variable $z_{ij}$
  (also to variables from $\mathfrak{k}$) the correspond the multiplicator  $\delta_i^{-1}\delta_j$. We suggest that $\delta_0=1$, $\delta_{1}=\delta_{-1}^{-1}$. The multiplicators are multiplied onto the function
    $\delta_{-n}^{m_{-n}}...\delta_{-1}^{m_{-1}}$, corresponding to the highest weight. The degree of $\delta_{-1}$ is the  weight $\Delta(m)_{-1}$.

 From the structure of $f$ one sees that the change of transformation of the weight $[m_{-n,n},...,m_{-1,n}]$ under the action of  $f$ equals to the sum of transformations under the action of $c$ and $f_0$. The transformation under the action of  $c$ was investigated when we considered the case $g_2$, it equals

 \begin{align*}
&-2\sum_{i=-2}^{T}m'_{i,k}+\sum_{i=-2}^{-1}m_{i,n}+m_{-2,n-1} \text { where }T=-1\text{ or }-2 \text{ in the case } \mathfrak{sp}_{2n},\mathfrak{o}_{2n},\\&
-2\sum_{i=-2}^{-1}m'_{i,k}+\sum_{i=-2}^{-1}m_{i,n}+m_{-2,n-1}-\sigma_{-n}\text { in the case } \mathfrak{o}_{2n+1}
 \end{align*}

Since the transformation of the weight under th action of $f_0$ is the same for all  $g_n$ and  $\mathfrak{gl}_{n+1}$ then using result for  $\mathfrak{gl}_{n+1}$ we obtain that under the action of $f_0$ to $\Delta(\Gamma)_{-1}$ the following value is added

 \begin{align*}
-2\sum_{i=-n}^{-3}m'_{i,k}+\sum_{i=-n}^{-3}m_{i}+\sum_{i=-n}^{-3}m_{-2,n-1}
 \end{align*}

Adding the transformations of the weight corresponding to $c$ and $f_0$ we prove the Proposition.

 \endproof

\section{The proof of Proposition \ref{pp}}
\label{redel}

\begin{prop*}
The following equality takes place
\begin{align}
<(\bar{m})_{max} \mid F_{\pm 1,-2} \mid (m)_{max}
>=<(\bar{m})_{red} \mid F_{\pm 1,-2} \mid (m)_{red}>_{red}
\end{align}
\end{prop*}

The equality is proved using the Wigner-Eckart theorem. One has
\begin{align}
\begin{split}&
<(\bar{m})_{max} \mid F_{-1,-2} \mid
(m)_{red}>_{max}=<(\bar{m})_{red}\mid F_{-1,-2} \mid
(m)_{red}>\cdot\\&\cdot<\begin{pmatrix}[\bar{m}]_{n-1}\\max
\end{pmatrix}\mid \begin{pmatrix}j \\ [1\dot{0}]_{n-1}\\-2\end{pmatrix} \mid \begin{pmatrix}(m)_{red}\\max
\end{pmatrix}>,
\end{split}
\end{align}

where the Wigner coefficient equals to $1$ if there exists $j$,
such that $[\bar{m}]_{n-1}=[m]_{n-1}+[0,...,1_{\text{ at the place
}j},...,0]$, and equals zero otherwise.  Thus one has
\begin{align} <(\bar{m})_{red} \mid F_{-1,-2} \mid
(m)_{red}>_{red}=<(\bar{m})_{max}\mid F_{-1,-2} \mid (m)_{max}>.
\end{align}

\endproof

\section{The proof of Proposition \ref{predl}}
\label{docpredl}

Let us be given a diagram $(m)$ that define a vector in s  $g_{n-1}$-representation
  with the highest weight $[m_{-n,n-1},...,m_{-2,n-1}]$.

A polynomial on the group $Z_{n-1}$ corresponds to this vector. Consider it as a polynomial on a bigger group  $Z_n$.

\begin{prop*}
 $g_{n}$-representation with the highest weight $[m_{-n,n-1},...,m_{-2,n-1},0]$. The corresponding tableau is of
 type $\begin{pmatrix}max\\ m\end{pmatrix}$, this is a  $g_{n}$-tableau that is obtained form
    $(m)$ by adding two maximum row.
\end{prop*}

\proof One has $$(m)=\zeta_1...\zeta_t(max),$$ where $\zeta_i\in
Z^{-}_{n-1}$ and $(max)$ is the highest vector.

Let use the realization on the group $Z$. The   highest vector is the function that equals identically to one.

One has $(m)=(
T_{\zeta_1}...T_{\zeta_t}1)(z)=\alpha_{n-1}(\widetilde{z\zeta_1...\zeta_t})$.

Also one has

\begin{align}
\begin{pmatrix} max\\ m \end{pmatrix}=\zeta_1...\zeta_t(max),\end{align} where $\zeta_i\in Z^{-}_{n-1}\subset Z^{-}_n$.

 In the space of functions on   $Z_n$, one has

\begin{align}
\begin{pmatrix} max\\ m \end{pmatrix}=(T_{\zeta_1}...T_{\zeta_t} 1)(z)=\alpha_{n}(\widetilde{z\zeta_1...\zeta_t}   )
\end{align}

If $\alpha_n$ corresponds to the highest weight
$[m_{-n,n-1},m_{-n+1,n-1},...,m_{-2,n-1},0]$, then
 $\alpha_{n-1}(\widetilde{z\zeta_1...\zeta_t})=\alpha_{n}(\widetilde{z\zeta_1...\zeta_t}
)$.

Thus the polynomials corresponding to  $(m)$ and $\begin{pmatrix} max\\ m
\end{pmatrix}$ coincide.
\endproof

\end{document}